\documentclass[12pt]{amsart}
\usepackage{amscd,amssymb,latexsym}
\usepackage{fullpage}

\DeclareMathOperator*{\hocolimit}{\mathrm{holim}}
\newcommand{\hcl}[1]{{\displaystyle
            \hocolimit_{\substack{\hbox to 20pt{\rightarrowfill} \\ #1}}\,}}

\newcommand{\growth}{\mathrm{growth}}

\newcommand{\fg}{\mathcal{FG}}

\newcommand{\HOmega}{H^{\Omega}}
\newcommand{\omegabar}{\overline{\omega}}

\newcommand{\ext}{\mathrm{Ext}}     
\newcommand{\tor}{\text{Tor}}       
\newcommand{\ba}{\boxtimes}         
\newcommand{\chains}{\mathrm{C}}    
\newcommand{\op}{\mathrm{op}}       
\newcommand{\cell}{\mathrm{Cell}}   
\newcommand{\po}{\mathcal{PO}}      
\newcommand{\co}{\mathcal{CO}}      

\newcommand{\shr}{\chi}             
\newcommand{\cou}{c}                

\newcommand{\Derived}{Ho}   

\newcommand{\A}{\mathcal{A}}	    %
\newcommand{\ee}{\mathcal{E}}       %
\newcommand{\cc}{\mathcal{C}}       %

\newcommand{\phat}{_p^{\wedge}}

\newcommand{\Mdef}[2]{\newcommand{#1}{\relax \ifmmode #2 \else $#2$\fi}}

\newcommand{\finbuilds}{\models}
\newcommand{\finbuiltby}{=\!\!|}
\newcommand{\builds}{\vdash}

\newcommand{\tensor}{\otimes}

\newcommand{\map}{\mathrm{map}}
\newcommand{\Hom}{\mathrm{Hom}}

\newcommand{\Ext}{\mathrm{Ext}}
\Mdef{\bhom}{\mathbf{\hat{H}om}}

\Mdef{\Mod}{\mathrm{mod}}

\newcommand{\st}{\; | \;}

\renewcommand{\thesubsection}{\thesection.\Alph{subsection}}

\newtheorem{thm}{Theorem}[section]
\newtheorem{lemma}[thm]{Lemma}
\newtheorem{prop}[thm]{Proposition}
\newtheorem{cor}[thm]{Corollary}

\theoremstyle{definition}

\newtheorem{defn}[thm]{Definition}

\newtheorem{example}[thm]{Example}

\newtheorem{remark}[thm]{Remark}
\newtheorem{problem}[thm]{Problem}
\newtheorem{conj}[thm]{Conjecture}

\newcommand{\qqed}{\qed \\[1ex]}
\renewenvironment{proof}[1][\hspace*{-.8ex}]{\noindent {\bf Proof #1:\;}}{\qqed}

\newcommand{\cosusp}[1]{\Sigma_{#1}}

\Mdef{\PH} {\Phi^H}
\Mdef{\PK} {\Phi^K}
\Mdef{\PL} {\Phi^L}
\Mdef{\PT} {\Phi^{\T}}

\Mdef{\ef}{E{\cF}_+}
\Mdef{\etf}{\widetilde{E}{\cF}}
\Mdef{\eg}{E{G}_+}
\Mdef{\etg}{\tilde{E}{G}}

\Mdef{\infl}{\mathrm{inf}}
\Mdef{\defl}{\mathrm{def}}
\Mdef{\res}{\mathrm{res}}
\Mdef{\ind}{\mathrm{ind}}
\Mdef{\coind}{\mathrm{coind}}

\Mdef{\univ}{\mathcal{U}}

\Mdef{\Fp}{\mathbb{F}_p}
\Mdef{\Zpinfty}{\Z /p^{\infty}}
\Mdef{\Zpadic}{\Z_p^{\wedge}}

\newcommand{\bi}{\begin{itemize}}
\newcommand{\be}{\begin{enumerate}}
\newcommand{\bc}{\begin{center}}
\newcommand{\bd}{\begin{description}}
\newcommand{\ei}{\end{itemize}}
\newcommand{\ec}{\end{center}}
\newcommand{\ed}{\end{description}}

\newcommand{\lra}{\longrightarrow}
\newcommand{\lla}{\longleftarrow}
\newcommand{\stacklra}[1]{\stackrel{#1}\lra}

\Mdef{\we}{\mathbf{we}}
\Mdef{\fib}{\mathbf{fib}}
\Mdef{\cof}{\mathbf{cof}}
\Mdef{\BI}{\mathcal{BI}}

\newcommand{\ann}{\mathrm{ann}}

\Mdef{\B}{\mathbb{B}}
\Mdef{\C}{\mathbb{C}}
\Mdef{\D}{\mathbb{D}}
\Mdef{\E}{\mathbb{E}}
\Mdef{\T}{\mathbb{T}}
\Mdef{\F}{\mathbb{F}}
\Mdef{\G}{\mathbb{G}}
\Mdef{\I}{\mathbb{I}}
\Mdef{\N}{\mathbb{N}}
\Mdef{\Q}{\mathbb{Q}}
\Mdef{\R}{\mathbb{R}}
\Mdef{\bbS}{\mathbb{S}}
\Mdef{\Z}{\mathbb{Z}}

\Mdef{\bA}{\mathbb{A}}
\Mdef{\bB}{\mathbb{B}}
\Mdef{\bC}{\mathbb{C}}
\Mdef{\bD}{\mathbb{D}}
\Mdef{\bE}{\mathbb{E}}
\Mdef{\bF}{\mathbb{F}}
\Mdef{\bG}{\mathbb{G}}
\Mdef{\bH}{\mathbb{H}}
\Mdef{\bI}{\mathbb{I}}
\Mdef{\bJ}{\mathbb{J}}
\Mdef{\bK}{\mathbb{K}}
\Mdef{\bL}{\mathbb{L}}
\Mdef{\bM}{\mathbb{M}}
\Mdef{\bN}{\mathbb{N}}
\Mdef{\bO}{\mathbb{O}}
\Mdef{\bP}{\mathbb{P}}
\Mdef{\bQ}{\mathbb{Q}}
\Mdef{\bR}{\mathbb{R}}
\Mdef{\bS}{\mathbb{S}}
\Mdef{\bT}{\mathbb{T}}
\Mdef{\bU}{\mathbb{U}}
\Mdef{\bV}{\mathbb{V}}
\Mdef{\bW}{\mathbb{W}}
\Mdef{\bX}{\mathbb{X}}
\Mdef{\bY}{\mathbb{Y}}
\Mdef{\bZ}{\mathbb{Z}}

\Mdef{\cA}{\mathcal{A}}
\Mdef{\cB}{\mathcal{B}}
\Mdef{\cC}{\mathcal{C}}
\Mdef{\mcD}{\mathcal{D}} 
\Mdef{\cE}{\mathcal{E}}
\Mdef{\cF}{\mathcal{F}}
\Mdef{\cG}{\mathcal{G}}
\Mdef{\mcH}{\mathcal{H}} 
\Mdef{\cI}{\mathcal{I}}
\Mdef{\cJ}{\mathcal{J}}
\Mdef{\cK}{\mathcal{K}}
\Mdef{\mcL}{\mathcal{L}}

\Mdef{\cM}{\mathcal{M}}
\Mdef{\cN}{\mathcal{N}}
\Mdef{\cO}{\mathcal{O}}
\Mdef{\cP}{\mathcal{P}}
\Mdef{\cQ}{\mathcal{Q}}
\Mdef{\mcR}{\mathcal{R}}
\Mdef{\cS}{\mathcal{S}}
\Mdef{\cT}{\mathcal{T}}
\Mdef{\cU}{\mathcal{U}}
\Mdef{\cV}{\mathcal{V}}
\Mdef{\cW}{\mathcal{W}}
\Mdef{\cX}{\mathcal{X}}
\Mdef{\cY}{\mathcal{Y}}
\Mdef{\cZ}{\mathcal{Z}}

\Mdef{\At}{\tilde{A}}
\Mdef{\Bt}{\tilde{B}}
\Mdef{\Ct}{\tilde{C}}
\Mdef{\Et}{\tilde{E}}
\Mdef{\Ht}{\tilde{H}}
\Mdef{\Kt}{\tilde{K}}
\Mdef{\Lt}{\tilde{L}}
\Mdef{\Mt}{\tilde{M}}
\Mdef{\Nt}{\tilde{N}}
\Mdef{\Pt}{\tilde{P}}

\Mdef{\tA}{\tilde{A}}
\Mdef{\tB}{\tilde{B}}
\Mdef{\tC}{\tilde{C}}
\Mdef{\tE}{\tilde{E}}
\Mdef{\tH}{\tilde{H}}
\Mdef{\tK}{\tilde{K}}
\Mdef{\tL}{\tilde{L}}
\Mdef{\tM}{\tilde{M}}
\Mdef{\tN}{\tilde{N}}
\Mdef{\tP}{\tilde{P}}

\Mdef{\ft}{\tilde{f}}
\Mdef{\xt}{\tilde{x}}
\Mdef{\yt}{\tilde{y}}

\Mdef{\Ab}{\overline{A}}
\Mdef{\Bb}{\overline{B}}
\Mdef{\Cb}{\overline{C}}
\Mdef{\Db}{\overline{D}}
\Mdef{\Eb}{\overline{E}}
\Mdef{\Fb}{\overline{F}}
\Mdef{\Gb}{\overline{G}}
\Mdef{\Hb}{\overline{H}}
\Mdef{\Ib}{\overline{I}}
\Mdef{\Jb}{\overline{J}}
\Mdef{\Kb}{\overline{K}}
\Mdef{\Lb}{\overline{L}}
\Mdef{\Mb}{\overline{M}}
\Mdef{\Nb}{\overline{N}}
\Mdef{\Ob}{\overline{O}}
\Mdef{\Pb}{\overline{P}}
\Mdef{\Qb}{\overline{Q}}
\Mdef{\Rb}{\overline{R}}
\Mdef{\Sb}{\overline{S}}
\Mdef{\Tb}{\overline{T}}
\Mdef{\Ub}{\overline{U}}
\Mdef{\Vb}{\overline{V}}
\Mdef{\Wb}{\overline{W}}
\Mdef{\Xb}{\overline{X}}
\Mdef{\Yb}{\overline{Y}}
\Mdef{\Zb}{\overline{Z}}

\Mdef{\db}{\overline{d}}
\Mdef{\hb}{\overline{h}}
\Mdef{\qb}{\overline{q}}
\Mdef{\rb}{\overline{r}}
\Mdef{\tb}{\overline{t}}
\Mdef{\ub}{\overline{u}}
\Mdef{\vb}{\overline{v}}

\Mdef{\hc}{\hat{c}}
\Mdef{\he}{\hat{e}}
\Mdef{\hf}{\hat{f}}
\Mdef{\hA}{\hat{A}}
\Mdef{\hH}{\hat{H}}
\Mdef{\hJ}{\hat{J}}
\Mdef{\hM}{\hat{M}}
\Mdef{\hP}{\hat{P}}
\Mdef{\hQ}{\hat{Q}}

\Mdef{\thetab}{\overline{\theta}}
\Mdef{\phib}{\overline{\phi}}

\Mdef{\uA}{\underline{A}}
\Mdef{\uB}{\underline{B}}
\Mdef{\uC}{\underline{C}}
\Mdef{\uD}{\underline{D}}

\Mdef{\bolda}{\mathbf{a}}
\Mdef{\boldb}{\mathbf{b}}
\Mdef{\boldD}{\mathbf{D}}

\Mdef{\fm}{\frak{m}}
\Mdef{\fp}{\frak{p}}

\Mdef{\eps}{\epsilon}

\input{xypic}

\begin{document}
\title{Complete intersections and mod $p$ cochains.}

\author{D.J.Benson}
\address{Department of Mathematics, University of Aberdeen, Aberdeen AB24 3UE,
UK}
\email{bensondj@maths.abdn.ac.uk}

\author{J.P.C.Greenlees}
\address{School of Mathematics and Statistics, Hicks Building,
Sheffield S3 7RH, UK}
\email{j.greenlees@sheffield.ac.uk}
\author{S.Shamir}
\address{Department of Mathematics, University of Bergen, 5008 Bergen, Norway}
\email{shoham.shamir@math.uib.no}
\date{}

\begin{abstract}
We give homotopy invariant definitions corresponding to  three well
known properties of complete intersections, for the ring, the module
theory and the endomorphisms of the residue field, and we investigate
them for the mod $p$ cochains on a space, showing that suitable
versions of the second and third are equivalent and that the first is
stronger. We are particularly interested in classifying spaces of
groups, and we give a number of examples. The case of rational
homotopy  theory is treated in \cite{qzci}, and there are some
interesting  contrasts. 
\end{abstract}

\thanks{The research was partially supported by EPSRC Grant number
EP/E012957/1. }
\maketitle

\tableofcontents



\section{Introduction}
\label{sec:Intro}
\subsection{The context.}
In algebraic geometry, the best behaved varieties are subvarieties
of an affine space which are specified by the right number of
equations:
if they are of codimension $c$ then only $c$ equations are required.
Considering this locally, we may work in commutative algebra.
A commutative local ring $R$ is a {\em complete intersection (ci)}
if its completion
is the quotient of a regular local ring $Q$ by a regular sequence,
$f_1,f_2, \ldots, f_c$. We will suppose that $R$ is complete, so that
$$R=Q/(f_1,\ldots, f_c).$$
The smallest possible value of $c$ (as $Q$ and the regular sequence vary)
is called the {\em codimension} of $R$.

If $R$ is a commutative Noetherian local ring with residue field $k$
that  is ci of codimension $c$, one may construct
a resolution of any finitely generated module growing like a polynomial
of degree $c-1$. In particular the ring $\Ext_R^*(k,k)$ has polynomial growth
(we say that $R$ is {\em gci}).
Perhaps the most striking result about ci rings is the theorem of Gulliksen
\cite{Gulliksenzci}
which states that this characterises ci rings so that the
ci and gci conditions are equivalent for local rings.

In fact one may go further and show that the resolutions are constructed
in an eventually multi-periodic fashion. In particular,
for a ci local ring $\Ext_R^*(k,k)$ is finite as a module over a
commutative Noetherian ring \cite{Gulliksen2}.
This in turn opens the way to
the theory of support varieties for modules over a ci ring \cite{AB}.

\subsection{The aspiration.}
We are interested in extending the notion of ci rings to commutative
differential graded algebras (DGAs) and commutative ring spectra. Indeed
we have a particular interest in studying the cochains $C^*(BG;k)$ on the
classifying space of a finite group $G$, where $k$ is a field of characteristic
$p$, partly because of the consequences for the cohomology ring $H^*(BG;k)$.
We would like to follow the model of \cite{DGI1},
which considers the Gorenstein condition. In fact, it shows $C^*(BG;k)$
is Gorenstein in a homotopical sense for all finite groups $G$. This structural result
then establishes the existence of a local cohomology theorem for
$H^*(BG;k)$. An immediate corollary is the result of Benson and Carlson
\cite{BC1}
that if $H^*(BG;k)$ is Cohen-Macaulay it is
also Gorenstein, and the fact that in any case $H^*(BG;k)$ is Gorenstein
in codimension 0 in the sense that its localization at any minimal
prime is Gorenstein.

By contrast, we only expect the complete intersection condition on
$C^*(BG)$ to hold for a small subclass of groups $G$, and we expect
the structural implications for $H^*(BG)$ to be at a more subtle
level.  To  explain this, we note that by the Eilenberg-Moore theorem,
the counterpart of the Ext algebra $\Ext_R^*(k,k)$ for $R=C^*(BG)$ is
the loop space homology of the $p$-completed classifying space,
$H_*(\Omega(BG_p^{\wedge}))$. Of course if $G$ is a $p$-group, this is
simply the group ring $kG$ in degree 0. More generally, it is known
to be of polynomial growth in certain cases (for instance $A_4$ or $M_{11}$
in characteristic 2) and
R.~Levi \cite{Levi1,Levi2,Levi3,Levi4}
has proved there is a dichotomy between small growth and large growth,
and given examples where the growth is exponential. Evidently groups
whose $p$-completed classifying spaces have
loop space homology that has exponential growth cannot be spherically resolvable,
so Levi's groups disproved a conjecture of F.Cohen.

We give homotopy invariant versions of all three characterizations of
the ci condition:
\begin{description}
\item[(sci)] the `regular ring modulo
regular sequence' condition,
\item[(mci)] the `modules have eventually
multiperiodic resolutions' and
\item[(gci)] polynomial growth of the Ext algebra
\end{description}
Before we do so, we need to give a counterpart of the Noetherian
condition, which in effect corresponds to Noether normalization. 
In crude terms, ignoring the variety of variants, 
we show that under the  Noetherian condition the mci and gci
conditions are equivalent, and both are implied by the strictly
stronger sci condition. These definitions acquire their interest
because of the variety of examples available, and the insights that
our results give. 

 We note that the Avramov-Quillen
characterization of ci rings in terms of Andr\'e-Quillen homology
does not work for cochains on a space in the mod $p$ context since Mandell has shown \cite{Mandell}
that the topological  Andr\'e-Quillen cohomology
vanishes rather generally in this case.

\subsection{Relation to other papers in the series.}
In \cite{cazci} two of us observed that the work of Eisenbud
\cite{Eisenbudzci} and Avramov-Buchweitz \cite{AB} allowed us to 
formulate a homotopy invariant version of the mci condition  (the zci
condition). The work of Gulliksen \cite{Gulliksenzci} 
shows that it is  equivalent to the sci and gci conditions for classical commutative local rings. The
lectures of Avramov \cite{AvramovCRM} give a very accessible account of 
the relevant material making it clear that the zci condition is
equivalent to the sci and gci conditions for classical commutative
local rings.

In \cite{qzci} two of us together with K.Hess formulated
a homotopy invariant version of the sci condition and applied it
in the rational homotopy theory of simply connected spaces, by taking
$R=C^*(X;\Q)$ for some strictly commutative model of the rational
cochains. The Noetherian condition is then simply that $H^*(X;\Q)$ is
Noetherian.  The sci condition implies the gci condition, and under
a strong Noetherian condition sci  is equivalent to gci. However, even
when they are equivalent, the zci condition is strictly
stronger. Accordingly, we formulated a weakened homotopy invariant
condition (eci)  which is another counterpart of the mci condition and
coincides with the zci condition for ungraded commutative rings. We then
showed that for strongly Noetherian rational DGAs the sci, eci and gci
conditions are equivalent.  For finite complexes, the three conditions
are closely related to the classical notion of an elliptic space.

In the present paper we consider the same conditions as in the
rational case, and study the relationship between them when applied
to $R=C^*(X;k)$ for a field $k$ of characteristic $p$, and a connected
$p$-complete space $X$ (with particular interest in the case
$X=BG_p^{\wedge}$). At the superficial level there is a technical
difficulty in that there is no commutative DGA model for $C^*(X;k)$
in general, but this difficulty is easily circumvented by the use of
commutative ring spectra, at the cost of having to work in a more
sophisticated technical setting. There is also the difficulty that
in mod $p$ homotopy theory one cannot expect to make complete
algebraic calculations. Neither of these are essential
differences. On the other hand, a really significant difference comes
from the fact that in the rational case the centre of $H_*(\Omega X)$
can have a smaller growth rate than $H_*(\Omega X)$ itself. By
contrast,  in the mod $p$ context work of F\'elix-Halperin-Thomas
\cite{FHT} shows that this does not happen. We are therefore able
to prove results which are in some ways a little stronger, and it
seems possible that the eci and zci conditions are again equivalent.


\subsection{Organization of the paper.}

In Section \ref{sec:introcochains} we introduce our philosophy in studying the cochain
algebras: the idea is to formulate conditions in commutative algebra so
that they are homotopy invariant  and then to say that a
space $X$ has Property P if the commutative ring spectrum $C^*(X;k)$
has Property P. This picks out interesting classes of spaces.

In Section \ref{sec:regular}  we record results about regular spaces, introducing
some terminology that becomes essential in discussions of complete
intersections.  There are two further sections of prerequisites.
In Section \ref{sec:noethfg}   we introduce two finiteness conditions that we need,
analogous to the Noetherian condition on rings and finite generation
for modules. In Section \ref{sec:bimod} we explain how bimodules give
rise to endomorphisms of the module category and how they are related
to  natural constructions on the module category.

This equips us to define in Section \ref{sec:cispaces}
a number of variants of the ci condition
for spaces: the sci condition is a structural condition like
the basic definition of a ci local ring, and it is the strongest.
The eci condition is analogous to having a multiperiodic module theory
and the gci condition is a growth condition on the loop space homology.

In the Sections \ref{sec:shyperiszhyper} to \ref{sec:gciiseci}  we explain how these
definitions are related.  Very roughly speaking Section
\ref{sec:shyperiszhyper} explains how sci implies eci,
Section \ref{sec:eciisgci} shows that eci implies gci,  and
Section \ref{sec:gciisGor} shows that gci implies Gorenstein.
Finally, the most surprising fact is that the finiteness condition gci
does have strong structural implications for the module categories.
Section \ref{sec:shearing} explains some basic properties of maps
between Hochschild cohomology groups, in preparation for
Sections \ref{sec:gciiszci}  and \ref{sec:gciiseci}
which show that gci implies eci (and also zci if the space is
finite).

Finally, we finish with a section of examples coming from
representation theory, and pose some questions.

We are grateful to N.Castellana, W.G.Dwyer and A.Gonzalez for a number
of valuable conversations related to these ideas. 

\subsection{Grading conventions.}
We will have cause to discuss homological and cohomological gradings.
Our experience is that this is a frequent source of confusion, so we adopt
the following conventions. First, we refer to lower gradings as {\em degrees}
and upper gradings as {\em codegrees}. As usual, one may convert gradings to
cogradings via the rule $M_n=M^{-n}$.
Thus both chain complexes and cochain
complexes have differentials of degree $-1$ (which is to say, of codegree $+1$).
This much is standard.  However, since we need to deal with both homology
and cohomology it is essential to have separate notation for
homological suspensions ($\Sigma^i$) and  cohomological suspensions
($\Sigma_i$): these are defined by
$$(\Sigma^iM)_n=M_{n-i} \mbox{ and } (\Sigma_iM)^n=M^{n-i}.$$
Thus, for example, with reduced chains and cochains  of a based
space $X$,  we have
$$C_*(\Sigma^iX)=\Sigma^iC_*(X) \mbox{ and } C^*(\Sigma^iX)=\Sigma_iC^*(X).$$


\section{Commutative algebra for spaces}
\label{sec:introcochains}

\subsection{Philosophy.}
For a space $X$ we want to work with a model of the cochains $R=C^*(X;k)$
which behaves like a commutative ring in the sense that its category of
modules has a symmetric monoidal tensor product $\tensor_R$, and so that
we can form the derived category of $R$-modules. We then use algebraic
behaviour of this commutative ring to pick out interesting classes of
spaces. In accordance with the principle that $C^*(X;k)$ is a sort of ring
of functions on $X$, we simplify terminology and say that $X$ has a property
P if the commutative ring $C^*(X;k)$ has the property P.

\subsection{Building}
We will be working in various triangulated categories. We say that $A$
{\em finitely builds} $B$ (or that $B$ is {\em finitely built by} $A$) if $B$ may be
formed from $A$ by finitely many cofibre sequences and retracts
(in other words, $B$ is in the thick subcategory generated by $A$). We
then write $A\finbuilds B$. 

Similarly if arbitrary coproducts are also permitted we say that $A$
{\em builds} $B$ (in other words, $B$ is in the localizing subcategory
generated by $A$). In this case we  write $A\builds B$

When we work in a category of $R$-modules, we say $M$ is {\em small} if 
it is finitely built by $R$ ($R\finbuilds M$). We say that $M$ is {\em
virtually-small} \cite{DGI2} if  it finitely builds a non-trivial small object
$W$, called the {\em witness}: $M\finbuilds W$, $R\finbuilds W\not \simeq 0$. We say that $M$ is 
{\em proxy-small} \cite{DGI1} if  it is virtually small and some
witness $W$ builds $M$.

\subsection{Convenient models.}
If $k$ is a field of characteristic zero, the simplicial de Rham
complex gives a strictly commutative model for the cochains, and the
commutative algebra of this DGA has been extensively investigated by rational
homotopy theorists \cite{FHT}, 
and from the present point of view in \cite{qzci}.
On the other hand,  it is
well known that in general Steenrod operations give an obstruction to a
natural commutative DGA model.
Fortunately \cite{EKMM,HSS}, 
there is a commutative ring
spectrum model of functions from $X$ to the Eilenberg-MacLane spectrum $Hk$,
so we use the notation
$$C^*(X;k) := \map (X,Hk) .$$
The category of module spectra over this ring spectrum has a model structure
with weak equivalences given by homotopy isomorphisms, and we write
$Ho(C^*(X;k))$ for the homotopy category: this is a tensor triangulated
category that provides a suitable setting for investigating the commutative
algebra for spaces.

This terminology is consistent in the sense that if we take $X$ to be
a point,  we have an equivalence
$$Ho(Hk) \simeq D(k), $$
of tensor triangulated categories, between the topological and
algebraic derived categories
\cite{Shipley}. 
If $k$ is replaced by a non-commutative ring
there is a  similar equivalence of triangulated categories.

We will not be using special properties of the models, so it is not necessary
to give further details, but \cite{spectra} gives an introduction and
guide to the available sets of foundations.

\subsection{Some analogies.}

At the most basic level, cofibre sequences
$$X \lra Y \lra Z$$
of pointed spaces induce (additive) exact sequences
$$C^*(X; k ) \lla C^*(Y;k) \lla C^*(Z; k) $$
of reduced cochains. On the other hand, fibrations
$$F \lra E \lra B$$
of spaces induced (multiplicative) exact sequences
$$ C^*(F;k)\stackrel{EM}\simeq C^*(E;k) \tensor_{C^*(B;k)} k \lla C^*(E;k) \lla C^*(B;k)$$
provided an Eilenberg-Moore theorem (EM) holds.

More generally, a homotopy pullback square
$$\begin{array}{ccc}
Z\times_XY&\lra & Z\\
\downarrow&& \downarrow\\
Y&\lra & X\\
\end{array}$$
induces a homotopy pushout square
$$\begin{array}{ccc}
C^*(Z\times_XY;k)&\lla & C^*(Z;k)\\
\uparrow&& \uparrow\\
C^*(Y;k)&\lla & C^*(X;k)\\
\end{array}$$
in the sense that
$$C^*(Z\times_XY;k)\simeq C^*(Z;k)\tensor_{C^*(X;k)}C^*(Y;k)$$
if the conditions of the Eilenberg-Moore theorem are satisfied
(for example if $X$ is $1$-connected or connected and
$p$-complete with
$\pi_1(X)$ a $p$-group and $p^N=0$ on $k$ for some $N$ \cite{Dwyer}).


In view of the importance of the Ext algebra, one particular case
will be especially significant for us.
\begin{prop}
\label{prop:EM}
We have an equivalence
$$C_*(\Omega X;k) \simeq \Hom_{C^*(X;k)}(k,k).$$
provided either (i)  $X$ is simply connected or
(ii) $X$ is connected and $p$-complete,
$\pi_1(X)$ is a $p$-group and $p^N=0$ on $k$ for some $N$.
\end{prop}

Because of the importance of this condition we will always assume that our
space satisfies either condition (i) or (ii). In particular, note that
if $X$ is the $p$-completion of a space with a finite fundamental group
(such as $B\Gamma$ for a compact Lie group) it satisfies the second
condition with $k=\Fp$.

\subsection{Conventions.}
Throughout we will be working over a field $k$ of
characteristic $p$, and we take   $k=\Fp$ for definiteness.
We repeat that the space $X$ is assumed to be $p$-complete, connected and to have
fundamental group a finite $p$-group. 

We will adapt our language accordingly, so that a space $X$ is $k$-finite
if $H_*(X;k)$ is finite dimensional, and a sequence of spaces
$F\lra E\lra B$ is a $k$-fibration if its $k$-completion is a fibration
in the usual sense.

We will often omit notation for the coefficients $k$, and say `finite'
when `$k$-finite' is intended. This amounts to a global assumption of
working in the category of $k$-complete spaces. 

The cochains, homology and cohomology of a space will always be taken
with coefficients in $k$, which will be omitted from the notation. Taking the cochains of a space $X$ results
in a commutative $Hk$-algebra $C^*(X)$ in the terminology of
\cite{EKMM}, which we will sometimes refer to simply as a  $k$-algebra. We will
use $k$ also to denote $Hk$ whenever there is no cause for confusion.

In this model, the homology is given by the homotopy groups
according to the formula  $\pi_*(C^*(X)) = H^{*}(X)$. We
say that a $C^*(X)$-module  $M$  is \emph{bounded above} if its
homotopy groups $\pi_*(M)$ are bounded above.

\section{Regular spaces}
\label{sec:regular}

Regular rings and spaces are not the main subject of this paper,
but it is essential to deal with regular rings first, 
since they provide  the basis for subsequent study: not only
do they provide the basis for an essential finiteness condition, but 
also complete interesections are defined as regular quotients of
regular rings.

\subsection{Commutative algebra.}
In commutative algebra there are three styles for a definition of
a regular local ring: ideal theoretic, in terms of the growth of
the Ext algebra and a derived version. Although these are equivalent
for classical commutative rings, we distinguish the definitions
for comparison with other contexts.

\begin{defn}
(i) A local Noetherian ring $R$ is {\em s-regular} if the maximal ideal
is generated by a regular sequence.

(ii) A local Noetherian ring $R$ is {\em g-regular} if $\Ext_R^*(k,k)$ is
finite dimensional.

(iii) A local Noetherian ring $R$ is {\em m-regular} if every finitely generated
module is small in the derived category $D(R)$.
\end{defn}

It is not hard to see that g-regularity is equivalent to m-regularity or that
s-regularity implies g-regularity. Serre proved that g-regularity
implies s-regularity \cite{Serre}, so the three conditions are equivalent.

\subsection{g-regularity for spaces.}

Of the commutative algebra definitions, the only one with a straightforward
counterpart for $C^*(X)$ is g-regularity. In view of Proposition
\ref{prop:EM} it takes the following form.

\begin{defn}
A space $X$  is {\em g-regular} if $H_*(\Omega X;k)$ is finite
dimensional.
\end{defn}

\begin{remark} Breaking our conventions for a moment, we might say
  that a not-necessarily-complete space $X'$ is globally regular if $\Gamma' =\Omega X' $
is a finite complex, so that $X=B\Gamma$ is the classifying space
of a finite loop space. For instance,
if $\Gamma$ is any compact Lie group we obtain such a space $X$.

Reverting to the standard situation that $k=\Fp$ and $X$ is
$p$-complete and connected,  the g-regularity condition
that $H_*(\Omega X; \Fp)$ is finite is precisely the condition
$\Gamma =\Omega X$ is
a $p$-compact group in the sense of Dwyer and Wilkerson \cite{DW}
with
classifying space $X=B\Gamma$. For example if $\Gamma'$ is a compact
Lie group with component group a $p$-nilpotent group, the $p$-completion of $B\Gamma'$
is an example, although there are many examples not of this form.
\end{remark}

\subsection{Some small objects.}

To discuss possible definitions of m-regularity we need to identify
some small objects.

\begin{lemma}
\label{lem:finfibissmall}
For a map $f: Y \lra X$ between $p$-complete spaces with fundamental
groups finite $p$-groups, $H_*(F(f))$ is finite
dimensional if and only if $C^*(Y)$ is small as a $C^*(X)$-module.
\end{lemma}

\begin{proof}
Suppose first that $H_*(F(f))$ is finite dimensional. 
Note first that since $\pi_1(X)$ is a finite $p$-group, the only
simple module over  $H_0(\Omega X)=k\pi_1X$ is the trivial module $k$.
Accordingly, the hypothesis implies that $k$ finitely builds $ C_*(F(f))$ as a
$C_*(\Omega X)$-module. Applying
$\Hom_{C_*(\Omega X)}(\cdot , k)$, we deduce from the Rothenberg-Steenrod
equivalence that $C^*(X)$ finitely builds $C^*(Y)$. In symbols,
$$k \finbuilds_{C_*(\Omega X)}  C_*(F(f))$$
and hence
$$C^*(X) \simeq  \Hom_{C_*(\Omega X)}(k,k) \finbuilds
\Hom_{C_*(\Omega X)}(C_*(F(f)),k) \simeq C^*(Y) .$$

The reverse implication is similar, but using the Eilenberg-Moore equivalence.
\end{proof}



\subsection{What should finitely generated mean?}

The m-regularity condition in commutative algebra
states that finitely generated modules are small. It would be nice to
have a counterpart of this for spaces, but it is not clear what should
play the role of finitely generated modules.

Our solution to this problem is to treat smallness as a {\em definition} of
`finitely generated' for modules over $g$-regular rings. We then use
$g$-regular spaces to define the notion in general. We pause to 
show that over $g$-regular spaces, at least for modules $M$ arising from 
maps of spaces, this notion of finitely generated has the familiar
form that $H^*(M)$ is finitely generated over $H^*(X)$. This also
holds in the rational context \cite{qzci}. 

\begin{lemma} (Dwyer) If $Y$ is a space with a map $f:Y\lra X$ 
to a regular space $X$ so that that $H^*(Y)$ is finitely generated
over $H^*(X)$ then in fact the fibre $F(f)$ is small. 
\end{lemma}

\begin{proof}
First,  by the fundamental work of Dwyer-Wilkerson on $p$-compact
groups \cite{DW}  $H^*(X)$ is Noetherian. Accordingly, since
$H^*(Y)$ is a finitely generated algebra over $H^*(X)$ it is also Noetherian. 
  
Finally, we argue that  $H^*(F)$ is finitely  generated over $H^*(Y)$. 
For this we consider the Serre spectral sequence of 
$$         \Omega X \lra  F \lra Y. $$
 The $E_2$-page has a finite number of rows, and each is a finitely  
generated  module over $H^*(Y)$.  It follows inductively that $ E_r$ 
has a finite number of rows,  each finitely generated over $H^*(Y)$. 
Since the spectral sequence collapses at a finite stage, this also
applies to $r=\infty$.  Piecing the rows together, we see that $H^*(F)$
is finitely generated  over $H^*(Y)$ as claimed.

Finally, we see that $H^*(F)$ is finite.  Indeed $H^*(Y)$ is Noetherian
and  the action of $H^*(Y)$ on
$H^*(F)$ factors  through the finite quotient ring
$k\tensor_{H^*X}H^*(Y)$. 
\end{proof}

\section{A Noetherian condition}
\label{sec:noethfg}

Before we turn to complete intersections there are two matters we need
to discuss. The first, which we deal with here, is the notion of a
normalizable  space, which is a counterpart of
the Noetherian condition.  This lets us use the ideas of the
previous section to give a notion of `finitely generated' module for
all normalizable spaces.  Thanks to some major theorems in the theory
of $p$-compact groups we can show that the simplest possible
definition has good properties.

\subsection{Normalizable spaces.}
In commutative algebra, it is natural to impose the Noetherian
finiteness assumption. One of the most useful consequences for $k$-algebras
is Noether normalization, stating that a Noetherian $k$-algebra is a
finitely generated module over a polynomial subring.

\begin{defn}
A space $X$ is {\em normalizable} if there is a connected g-regular space $B\Gamma$
and a map $\nu: X \lra B\Gamma$ so that the homotopy fibre
$F(\nu)$ is a $k$-finite complex. We say that $\nu$ is a
{\em normalization} and $F(\nu)$ is a {\em Noether fibre}.
The normalization is called \emph{polynomial} if $B\Gamma$ is
simply-connected and has polynomial cohomology.
\end{defn}


\begin{remark}
(i) 
By Lemma~\ref{lem:finfibissmall}, 
the condition that $F(\nu)$ is $k$-finite is equivalent to asking that $C^*(X)$
is small over $C^*(B\Gamma)$.

(ii) In the rational context, whenever $H^*(X)$ is Noetherian the space
$X$ is normalizable \cite{qzci}.

(iii) If $X$ is normalizable then $H^*(X)$ is Noetherian:
since $H^*(B\Gamma)$ is \cite{DW}
Noetherian the first remark implies $H^*(X)$ is a finitely generated
$H^*(B\Gamma)$-module.

(iv) Since our spaces $X$ are assumed connected, the Noether fibre of a
polynomial  normalization is connected. 
\end{remark}

\begin{example}
If $G$ is a compact Lie group (for example a finite group),
then $BG$ is polynomially normalizable by
choosing a faithful representation $G \lra U(n)$, giving a fibration
$U(n)/G \lra BG \lra BU(n)$.
\end{example}

It is clear that requiring the existence of a map $X \lra B\Gamma$
with $C^*(X)$ small over $C^*(B\Gamma)$ is a  finiteness condition, but
we need to give examples to see how stringent it is. For example,
it is natural to assume $H^*(X)$ is Noetherian, but this does not 
guarantee normalizability.

\begin{example} {\em (Castellana)}  A space $X$ for which
  $H^*(X;\Z_p^{\wedge})$ has unbounded $p$-torsion is not
  normalizable.  
Indeed, if $X\lra B\Gamma$ is a normalization, 
 it is shown by Dwyer-Wilkerson \cite{DW} that
$H^*(X; \Z_p^{\wedge})$ is finitely generated over the Noetherian
ring  $H^*(B\Gamma;\Z_p^{\wedge})$ and hence Noetherian. 

Various classes of such spaces with $H^*(X;\F_p)$ Noetherian are
known. A simple example is the 3-connected cover of  $S^3$. Others
come up naturally when considering homotopical generalizations of 
classifying spaces of groups. For instance, some of the classifying spaces of rank 2 Kac-Moody groups
described by Kitchloo \cite{Kitchloo}  have this property, as do the 
 Aguad\'e-Broto-Notbohm \cite[5.5 and 5.6]{ABN} spaces $X_k(r)$. 
\end{example}

\subsection{Finitely generated modules.}
Next we need to define a notion of finitely generated
modules, and for this we restrict to normalizable spaces $X$.
We choose a normalization $\nu: X\lra B\Gamma$ and  then take
$$\fg_{\nu} :=\{ C^*(Y) \st Y\lra X \stackrel{\nu}\lra B\Gamma
\mbox{ makes } C^*(Y) \mbox{ small over } C^*(B\Gamma)\}, $$
and
$$\fg (X) :=\bigcup_{\mathrm{normalizations}\;\; \nu}\fg_{\nu}.$$


In any case, it is convenient to be able
to compare different normalizations. Thus if $\nu_1, \nu_2$ are
two normalizations we can form
$$\{\nu_1, \nu_2\} : X \lra B\Gamma_1 \times B\Gamma_2.  $$
We may compare $\{\nu_1, \nu_2\}$ to $\nu_1$ in the diagram
$$\diagram
X \dto^= \rto^(0.3){\{ \nu_1,\nu_2\} } &B\Gamma_1 \times B\Gamma_2\dto \\
X \rto^{\nu_1} &B\Gamma_1.
\enddiagram$$

\begin{lemma}
Given two normalizations, $\nu_1, \nu_2$, the map
$\{\nu_1, \nu_2\}$ is also a normalization, and we have
$$\fg_{\nu_1} \subseteq \fg_{ \{\nu_1,\nu_2\} }. $$
\end{lemma}

\begin{proof}
For the first statement, we need only take iterated fibres
to obtain a fibration
$$\Gamma_2 \lra F(\{\nu_1, \nu_2\} ) \lra F(\nu_1)$$
from which the statement follows. For the second statement
we repeat the proof with $X$ replaced by $Y$ and the normalizations
replaced by their composites with $Y \lra X$.
\end{proof}

We next show that the class $\fg$ is independent of the chosen
normalization when the following conjecture holds.

\begin{conj} {\em (Linear representation.)}
\label{LRC}
The conjecture is that every $p$-compact group has a faithful linear
representation. More precisely that if
$B\Gamma$ is regular there is a map $B\Gamma \lra
BSU(n)$, for some $n$, whose homotopy fibre is $\Fp$-finite.
We understand that the conjecture has in fact been proved, but we 
state it as a conjecture so that the dependence of our work on it is
made clear. 

For the prime 2 the classification of $2$-compact
groups \cite{twocompact, twocompactMone, 
twocompactMtwo} shows that any 2-compact group is the product of
classical groups and copies of $DI(4)$. Ziemianski \cite{Ziemianski}
has constructed a faithful linear representation of $BDI(4)$. 

For odd primes $p$ this is a consequence of the classification of $p$-compact
groups \cite[1.6]{pcompact} using verifications in some exotic cases
by Castellana \cite{Castellana1,Castellana2}. 
\end{conj}

Note the Linear  Representation Conjecture (LRC) implies that
if  $\nu:X \lra B\Gamma$ is a normalization, then we may compose with
a faithful linear representation to obtain another normalization
$X \lra BSU(n)$.

Thus for every normalizable space $X$ we can always assume without
loss of generality that the normalization is $\nu:X \lra BSU(n)$. In
particular we may assume that our normalization is polynomial.

\begin{lemma}
Assuming the LRC \ref{LRC}, if  $X$ is a normalizable space then $Y\in \fg (X)$
if and only if $H^*(Y)$ is finitely generated as a module over $H^*(X)$.
\end{lemma}
\begin{proof}
Suppose $Y\in \fg_{\nu}$,  so that there are maps $Y \to X
\stackrel{\nu}\lra  B\Gamma$ making $C^*(Y)$ into a small $C^*(B\Gamma)$-module.
Since $H^*(B\Gamma)$ is Noetherian, this implies $H^*(Y)$ is finitely generated
over $H^*(B\Gamma)$ and hence also over $H^*(X)$.

Conversely, suppose given a map $Y \lra X$ such that $H^*(Y)$ is
finitely generated as an $H^*(X)$-module. Since $X$ is normalizable,
the LRC \ref{LRC} shows we can choose a normalization
$\nu: X \lra BSU(n)$ for some $n$. Now $H^*(Y)$ is a finitely
generated module over $H^*(BSU(n))$, and the Eilenberg-Moore spectral
sequence shows that the homotopy fibre of $Y \lra BSU(n)$ is $k$-finite.
Hence $Y\in \fg_{\nu}$.
\end{proof}

Even more is true.

\begin{lemma}
Assuming the LRC \ref{LRC}, if $X$ has a polynomial normalization $\nu:X \lra B\Gamma$
then $\fg_{\nu} = \fg$.
\end{lemma}

The proof is left as an exercise to the reader.

In the following sections we will give several definitions where we
require a polynomial normalization, so that we can rely on this
lemma.

\section{Bimodules and natural endomorphisms of $R$-modules}
\label{sec:bimod}

The second topic we need to treat before coming to complete
intersections is that of bimodules.

\subsection{The centre of the derived category of $R$-modules}
If $R$ is a commutative Noetherian ring and $M$ is a finitely
generated  $R$-module with an eventually $n$-periodic resolution,
comparing the resolution with its $n$-fold shift gives a 
map $M\lra \Sigma^nM$ in the derived category whose mapping cone
is  a small $R$-module. If $R$ is a hypersurface ring, any finitely generated module
has an eventually 2-periodic resolution. In fact the construction of
these resolutions can be made in a very uniform way. The lesson learnt 
from commutative algebra is that if we want to use this property 
to characterize hypersurfaces we need to use this uniformity. 
In fact the uniform construction can be formulated as a natural
transformation $1\lra \Sigma^n 1$ of the identity functor with a
mapping cone which is small on finitely generated modules, and it turns
out that the existence of such a transformation does characterize
hypersurface rings. 

By definition the centre $ZD(R)$ is the graded ring of all such
natural transformations of the identity functor.
There are various ways of constructing elements of the centre, and
various natural ways to restrict the elements we consider.
 Some of these work better than others, and it is the purpose of this
section is to introduce these ideas.

\subsection{Bimodules}

We consider a map $Q\lra R$, where $Q$ is regular and $R$ is small
over $Q$. We may then
consider $R^e=R\tensor_QR$, and $R^e$-modules are $(R|Q)$-bimodules.
The Hochschild cohomology ring is defined by
$$HH^*(R|Q)=\Ext_{R^e}^*(R,R).$$

If $f:X\lra Y $ is a map of  $(R|Q)$-bimodules, for any $R$-module $M$
we obtain a map  $f\tensor 1 : X\tensor_R M\lra Y\tensor_R M$ of
(left) $R$-modules.

The simplest way for us to use this is that if we have isomorphisms
$X\cong R$ and  $Y\cong \Sigma^n R$ as $R$-bimodules, the map
$f\tensor 1: M\lra \Sigma^n M$  is natural in $M$ and therefore gives
an element of codegree $n$ in $ZD(R)$: we obtain a map of rings
$$HH^n(R) =\Hom_{R^e}(R,\Sigma^n R)\lra ZD(R)^n.$$

Continuing, if $X\finbuilds_{R^e} Y$ then $X\tensor_R M\finbuilds_R
Y\tensor_R M$. In particular, if $X=R$ builds a small $R^e$-module
$Y$ then
$$ M=R\tensor_R M\finbuilds_R Y\tensor_R M \finbuiltby
R^e\tensor_RM=R\tensor_Q M.$$
Thus if $M$ is finitely generated (i.e., small over $Q$), this shows
$M$ finitely builds a small $R$-module.

We could then restrict the maps permitted in showing that $X
\finbuilds_{R^e}Y$. We could restrict ourselves to using maps of
positive codegree coming from Hochschild cohomology, and
say $X\finbuilds_{hh} Y$, more generally we could permit any maps of positive codegree from the centre
$ZD(R^e)$ and say $X\finbuilds_z Y$,  or we could relax further and
require only that all the maps involved in building are endomorphisms
 of non-zero degree for some object and say $X\finbuilds_e Y$.



We shall have occasion to use Koszul constructions of modules.
Given $z_1,...,z_c \in ZD(R)$ and an $R$-module $M$ we define
$M/z_1$ to be the homotopy cofibre of $\Sigma^n M \stacklra{z_1} M$ and
$M/z_1/\cdots/z_c$ is defined inductively. Given elements
$x_1,...,x_n\in HH^*(R|Q)$, we have already observed that these define elements in $ZD(R)$
(because $R$ is a unit for $\otimes_R$), and we note here that they
also define elements in $ZD(R^e)$ for the same reason. Accordingly, 
for a bimodule $M$ we may construct $M/x_1/\cdots/x_c$ as a bimodule.

\begin{example}
In the topological context we have a map $X\lra B\Gamma$ and
take $R=C^*(X), Q=C^*(B\Gamma)$ and $R^e=C^*(X\times_{B\Gamma}X)$.
The associated Hochschild cohomology can be abbreviated
$$HH^*(X|B\Gamma)=\pi_*(\Hom_{C^*(X\times_{B\Gamma}X)}(C^*(X), C^*(X)).$$
\end{example}

\section{Complete intersection spaces}
\label{sec:cispaces}

\subsection{The definition in commutative algebra}

In commutative algebra there are three styles for a definition of
a complete intersection ring: structural, in terms of the growth of
the Ext algebra and module theoretic. See \cite{AvramovCRM, cazci} for a more
complete discussion.

\begin{defn}
(i) A local Noetherian ring $R$ is an {\em s-complete intersection (sci)}
ring if $R=Q/(f_1,f_2, \ldots , f_c)$ for
some regular ring $Q$ and some regular sequence $f_1, f_2, \ldots , f_c$.
The minimum such $c$ (over all $Q$ and regular sequences) is called
the {\em codimension} of $R$.

(ii) A local Noetherian ring $R$ is  {\em gci}
if $\Ext_R^*(k,k)$ has polynomial growth. The  {\em g-codimension}
of $R$ is one more than the degree of the growth.

(iii) A local Noetherian ring $R$ is  {\em zci} \cite{cazci} if
there are elements $z_1,z_2,\ldots z_c\in ZD(R)$ of non-zero degree
so that
$M/z_1/z_2/ \cdots /z_c$ is small for all finitely generated modules
$M$. The minimum such $c$ is called the  {\em z-codimension} of $R$.
Similarly $R$ is {\em hhci} if the elements $z_i$ can be chosen to
come from Hochschild cohomology.
\end{defn}

\begin{remark}
The following condition has proved less useful.
If $R$ is a commutative ring or CDGA, it is said to be a
{\em quasi-complete intersection (qci)} \cite{DGI2}
if every finitely generated object  is virtually small.
\end{remark}

\begin{remark}
If $Q$ is a regular local ring with a map $Q\lra R$ making $R$ into a
small $Q$-module we may consider a number of bimodule conditions.

We say $R$ is {\em bci, eci, zbci} or {\em  hhbci} if $R$ finitely builds,
e-builds, z-builds or hh-builds a nontrivial small $R^e$-module.
Evidently these are increasingly strong conditions. Similarly,
it is clear that zbci implies zci and hhbci implies hhci.
\end{remark}


\begin{thm} \cite{cazci}
For a local Noetherian ring the conditions sci, gci, hhbci, hhci,
zbci, zci and eci are all  equivalent, and the corresponding codimensions are equal.
 These conditions imply the bci and qci conditions.
\end{thm}

It is a result of Shamash \cite{Shamash} that if $R$ is ci of codimension
$c$, one may construct a resolution of any finitely generated module growing
like a polynomial of degree $c-1$. Furthermore, the resolution is
constructed in \cite{Eisenbudzci} as a module over a polynomial ring on
$c$ generators of degee $-2$ which shows directly that $R$ is hhbci.
(the form of the resolution due to  Avramov-Buchweitz \cite{AB} and 
described in \cite[Section 9]{AvramovCRM} makes this very clear). 
Considering the module $k$ shows that the ring $\Ext_R^*(k,k)$ has polynomial
{\bf g}rowth. Perhaps the most striking result about
ci rings is the theorem of Gulliksen
\cite{Gulliksenzci} which states that gci  implies sci for commutative
local rings so that the
ci and gci conditions are equivalent.

\begin{remark}
In commutative algebra, Avramov \cite{AvramovAQ} proved Quillen's conjectured characterization
of complete intersections by the fact that the Andr\'e-Quillen cohomology
is bounded. The natural counterpart of this is false in homotopy theory
since, by results of Dwyer and Mandell \cite{Mandell}, the topological Andr\'e-Quillen
cohomology vanishes much too generally in characteristic $p$.

When $k$ is of characteristic 0, the Andr\'e-Quillen cohomology of $C^*(X)$
gives the dual homotopy groups of $X$, so Avramov's characterization corresponds
to the gci condition.
\end{remark}

\subsection{Definitions for spaces.}
In view of the fact that regular elements correspond to
spherical fibrations, adapting the above definitions for spaces is straightforward.

\begin{defn}
(i) A space $X$ is {\em spherically ci (sci)} if it is formed from a connected
g-regular space $B\Gamma$
using a finite number of spherical fibrations. More precisely, we require that
there is a g-regular space $X_0=B\Gamma$ and
fibrations
$$S^{n_1}\lra X_1 \lra X_0=B\Gamma, S^{n_2}\lra X_2 \lra X_1, \ldots ,
S^{n_c}\lra X_c \lra X_{c-1}$$
with $X=X_c$. The least such $c$ is called the {\em s-codimension} of
$X$.

(i)' A space $X$ is {\em weakly spherically ci (wsci)} if there is a
g-regular space $B\Gamma$ and a fibration
$$F\lra X \lra B\Gamma$$
and $F=F_1$ is spherically resolvable in the sense that there are fibrations
$$F_2\lra F_1 \lra S^{n_1}, F_3\lra F_2 \lra S^{n_2}, \ldots ,
*\lra  F_c\lra S^{n_c}$$
 The least such $c$ is called the {\em ws-codimension} of $X$.

(ii) A space  $X$ is a {\em gci} space
if $H^*(X)$ is Noetherian and $H_*(\Omega X)$ has polynomial growth.
The {\em g-codimension} of $X$ is one more than the degree of growth.


(iii) A space $X$ is a {\em zci} space if $X$ is normalizable and
there are elements $z_1, z_2, \ldots , z_c\in Z D(C^*(X))$ of non-zero degree
so that $C^*(Y)/z_1/z_2/ \cdots / z_c$ is small for all
$C^*(Y) \in \fg (X)$.

(iii)' A space $X$ is a {\em hhci} space if there is a polynomial normalization
$\nu:X \lra B\Gamma$ and $R=C^*(X)$ finitely $hh$-builds
a small $R^e$-module, where $R^e=C^*(X\times_{B\Gamma} X)$. In particular,
assuming the LRC \ref{LRC}, there are elements $z_1, z_2, \ldots , z_c\in HH^*(X|B\Gamma)$ of
positive codegree so that $C^*(Y)/z_1/z_2/ \cdots / z_c$ is small for all
$C^*(Y)\in \fg (X)$.

(iii)'' A space $X$ is a {\em eci} space if $X$ has a polynomial normalization
as above and $R$ finitely $e$-builds a small $R^e$-module. In particular, there
are homotopy cofibration sequences of $R^e$-modules
$$M_0\lra \Sigma^{n_0}M_0 \lra M_1,
M_1\lra \Sigma^{n_1}M_1 \lra M_2 , \ldots ,
M_{c-1}\lra \Sigma^{n_{c-1}}M_{c-1} \lra M_c  $$
with $n_i \neq 0$ and $M_c$ small over $R^e$.

\end{defn}

\begin{remark}
Two further variants have proved to be less useful.

We say $X$ is  {\em bci} space if $X$ is normalizable and
$C^*(X)$ is virtually small as a $C^*(X\times_{B\Gamma}X)$-module for some
g-regular space $B\Gamma$ and map $X\lra B\Gamma$ with $C^*(X)$ small over $C^*(B\Gamma)$.

We say $X$ is  {\em qci} space if $X$ is normalizable and
each $C^*(Y) \in \fg (X)$ is virtually small.

There are also further variants, $\omega$sci and w$\omega$sci
 where we are permitted to  use loop spaces on spheres rather than spheres.
These conditions arose in Levi's work. They are evidently weakenings
of sci and wsci, but they still imply gci.
\end{remark}

It is easiest to verify the gci condition. The simplest example gives a good supply.
\begin{example}
\label{eg:homciisgci}
If $H^*(X)$ is a complete intersection then $\Ext_{H^*(X)}^{*,*}(k,k)$
has polynomial growth, so that if the Eilenberg-Moore spectral
sequence converges, $X$ is gci. 
\end{example}

For normalizable spaces $X$, working over $k=\Fp$
we will establish the implications
$$sci \stackrel{}\Rightarrow wsci \stackrel{} \Rightarrow eci
\stackrel{} \Leftarrow\Rightarrow gci . $$
We will also show that for $k$-finite spaces $hhci
\stackrel{} \Leftarrow\Rightarrow gci$.
It is shown in \cite{qzci} that in the rational case
$sci \Rightarrow eci \Rightarrow gci$, and that 
 under an additional finiteness hypothesis
the three conditions sci, eci and gci are equivalent. 

The implication, that sci implies wsci is straightforward.  Indeed,
if $X$ is sci,  the fibre of  the composite
$$X=X_{c}\lra X_{c-1}\lra \cdots \lra X_1 \lra B\Gamma, $$
is clearly an iterated spherical fibration.  

In Section \ref{sec:eciisgci} we explain that  a rather straightforward
calculation with Hilbert series shows that eci implies gci.
In  Section \ref{sec:shyperiszhyper} we give a direct construction to
show that an s-hypersurface is a z-hypersurface. This argument
can be iterated if all the spherical fibrations have odd dimensional
spheres, but we give a less direct general argument.

The main result is that (assuming normalizability) the gci finiteness
condition implies the structural condition eci: this occupies
Sections \ref{sec:gciisGor} to \ref{sec:gciiseci}.

\begin{problem}
Give an example of a gci space which is not wsci or show that no such
space exists. 
\end{problem}



\subsection{Hypersurface rings.}

A hypersurface is a complete intersection of codimension 1. The first
three definitions adapt to define s-hypersurfaces, g-hypersurfaces and
z-hypersurfaces. The notion of g-hypersurface (i.e., the dimension of
the groups $\Ext^i_R(k,k)$ is bounded) may be strengthened to the
notion of p-hypersurface where we require
that they are eventually periodic, given by multiplication with an element
of the ring. All four of these conditions are equivalent
by results of Avramov.

One possible formulation of b-hypersurface would be to require that
$R$ builds a small $R^e$-module in one step (or equivalently,
that $R$ is a
z-hypersurface but $z$ arises from $HH^*(R|Q)$). This is also equivalent
to the above definitions.

Finally, we may say that $R$ is a q-hypersurface if every finitely
generated module $M$ has a self map with non-trivial small mapping cone.

\subsection{Hypersurface spaces.}

All six of these conditions have obvious formulations for spaces.
A space $X$ is an {\em s-hypersurface} if there is a fibration
$$S^n \lra X \lra B\Gamma$$
where $X \lra B\Gamma$ is a polynomial normalization. 
It is a {\em g-hypersurface}
if the dimensions of  $H_i(\Omega X)$ are bounded, and a {\em p-hypersurace}
if they are eventually periodic given by multiplication by an element of the
ring.

The space $X$ is an {\em hh-hypersurface} 
if $C^*(X)$  builds a small
$C^*(X\times_{B\Gamma}X)$-module in one step using an hh-map.
An e-hypersurface is the same as a hh-hypersurface. It is a {\em
  z-hypersurface} if there is an element $z$ of $ZD(C^*(X))$ of non-zero
degree so that  $C^*(Y)/z$ is small for every module $C^*(Y)$ in $\fg
(X)$. Finally,  $X$ is a {\em q-hypersurface} if every
module $C^*(Y)$ in $\fg (X)$ has a self map with non-trivial small mapping cone.


The point about this condition is that there are a number of interesting
examples (see Section \ref{sec:groupegs} below). We will show in Section
\ref{sec:shyperiszhyper} that an s-hypersurface is a z-hypersurface.

\section{s-hypersurface spaces and z-hypersurface spaces}
\label{sec:shyperiszhyper}

In the algebraic setting the remarkable fact is that modules
over hypersurfaces have eventually periodic resolutions, and hence that
they are hhci of codimension 1. The purpose of this section is to prove
a similar result for spaces. The result will be proved in a form that
provides the inductive step for the general result on ci spaces, but
the iteration will only show it is an eci space (rather than zci or hhci).

\begin{thm}
\label{thm:shyperiszhyper}
If $X$ is an s-hypersurface space with fibre sphere of dimension $\geq 2$
then $X$ is a z-hypersurface space.
\end{thm}

\begin{remark}
The proof will show a similar result holds for the total space of a circle fibration over a
connected regular space, but the definition of a z-hypersurface would
need to be adapted along the Jacobson radical lines of \cite{cazci} to give a uniform
statement.
\end{remark}

\subsection{Split spherical fibrations.}
\label{subsec:ssf}
The key in algebra was to consider  bimodules, for which we consider the
(multiplicative) exact sequence
$$R \lra R^e \lra R^e \tensor_{R} k , $$
where the first map is a monomorphism split by the map $\mu$ along which $R$ acquires
its structure as an $R^e$-module structure.
This corresponds to the pullback fibration
$$X \lla X\times_{B\Gamma}X \lla S^n, $$
split by the diagonal
$$\Delta : X \lra X \times_{B\Gamma }X$$
along which the cochains on $X$ becomes a bimodule.
To simplify notation, we consider a more general situation: a fibration
$$B \lla E \lla S^n$$
with section $s: B \lra E$. The case of immediate interest  is $B=X$,
$E=X\times_{B\Gamma}X$, where a $C^*(E)$-module is a $C^*(X)$-bimodule.

Since $s$ is a section of $p$ there is a fibration
$$\Omega S^n \lra B \stackrel{s}\lra E. $$
This gives the required
input for the following theorem. The strength of the result is that the
cofibre sequences are of $C^*(E)$-modules.

\begin{thm}
\label{thm:loopspherebuild}
Suppose given a fibration $\Omega S^n \lra B \stackrel{s}\lra E$
with $n\geq 2$.

(i) If $n$ is odd,
then there is a cofibre sequence of $C^*(E)$-modules
$$\Sigma_{n-1} C^*(B) \lla C^*(B) \lla C^*(E). $$

(ii)  If $n$ is even,
then there are cofibre sequences of $C^*(E)$-modules
$$C \lla C^*(B) \lla C^*(E) $$
and
$$\Sigma_{2n-2} C^*(B) \lla C \lla \Sigma_{n-1}C^*(E). $$
In particular the fibre of the composite
$$C^*(B)\lra C \lra \Sigma_{2n-2}C^*(B)$$
is a small $C^*(E)$-module constructed with one cell in codegree
0 and one in codegree $n-1$.
\end{thm}

\begin{remark}
The case $n=1$ is slightly different since $\Omega S^1\simeq \Z$.
In this case, the map $B \lra E$ is an infinite cyclic covering up to
homotopy. Taking appropriate models, $B$ is a free $\Z$-space with
$E=B/\Z$ and there is a cofibre sequence
$$C^*(E)\lra C^*(B) \stackrel{1-z}\lra C^*(B), $$
where $z$ is a generator of $\Z$.
\end{remark}

\subsection{Strategy}
We will first prove the counterparts
in cohomology by looking at the Serre spectral sequence
of the fibration from Part (i) and then lift the conclusion to
the level of cochains.

Note that in either case we obtain a cofibre sequence
$$ K\lla C^*(B) \lla \Sigma_a C^*(B)$$
of $C^*(E)$-modules with $K$ small.

In the motivating example we see that if $X$ is an s-hypersurface,
it is an e-hypersurface with $n \geq 2$ and so it is a z-hypersurface as required.
If $n=1$, the condition on non-triviality for a z-hypersurface needs
to be adapted along the lines of \cite{cazci}.
Thus Theorem \ref{thm:shyperiszhyper} follows from Theorem
\ref{thm:loopspherebuild}.

To describe the strategy in more detail, we consider the Serre spectral sequence
of the fibration in Part (i).

If $n$ is odd,  $H^*(\Omega S^n)$ is a free divided power algebra on one
generator $\Phi$ of codegree $n-1$.
If $n$ is even,  $H^*(\Omega S^n)$ is the tensor product of an
exterior algebra on one generator $\Phi$ of codegree $n-1$ and a
free divided power algebra on one generator $\Psi$ of codegree $2n-2$.

In either case, $H^*(\Omega S^n)$
corresponds precisely to the natural algebraic resolution of $R=C^*(B)$ over $C^*(E)$.
Indeed, the Serre spectral
sequence of this fibration gives precisely such a resolution, with
the generators of $H^*(\Omega S^n)$ giving an $H^*(E)$-basis for the
$E_2$-term. One might imagine realizing the associated
filtration by codegree in $H^*(\Omega S^n)$. Using the equivalence,
$$E\simeq E\Omega B \times_{\Omega B} \Omega S^n$$
the filtration would correspond to the skeletal filtration of $\Omega S^n$ provided this was a filtration of  $\Omega B$-spaces. We could then take
$$E^k=E\Omega B \times_{\Omega B} (\Omega S^n)^{((n-1)k)}$$
to realize the filtration, but the terms would have to be
$\Omega B$-spaces so it seems unlikely that this can be realized. In any
case we can make do with a little less.



To start with we consider the map
$$ C^*(B) \stackrel{s^*}\lla C^*(E)$$
induced by the section, using it to give $C^*(B)$ the structure of a $C^*(E)$-module. Thus $s^*$ is a map of $C^*(E)$-modules, and we can take its
mapping cone $C$ in the category of $C^*(E)$-modules.

By looking
at the Serre spectral sequence of the original spherical fibration,
we see that $\pi_*(C)\cong \cosusp{n-1}H^*(B)$ as $H^*(B)$-modules, and hence
$C$ is equivalent to $\cosusp{n-1}C^*(B)\simeq C^*(\Sigma^{n-1} B)$
as a $C^*(B)$-module.  However, we need to consider $\pi_*(C)$
not just as an $H^*(B)$-module, but as an $H^*(E)$-module.

If $n$ is odd, we will show $\pi_*(C)$ is $\cosusp{n-1}H^*(B)$.
If $n$ is even this need not be true, but we may repeat the
construction,  once more to obtain a $C^*(E)$-module $C_1$
so that $\pi_*(C_1)\cong \cosusp{2n-2}H^*(B)$ as an $H^*(E)$-module.

In either case we have a $C^*(E)$-module $M$ so that $\pi_*(M) \cong H^*(B)$
as $H^*(E)$-modules. We will show in Subsection
\ref{subsec:liftingtocochains} that the algebraic simplicity of the
$H^*(E)$-module  $H^*(B)$ is such that we can lift the equivalence
to the level of cochains.

If $n$ is odd, we therefore have a triangle
$$ \cosusp{n-1} C^*(B) \lla C^*(B) \lla C^*(E)$$
of $C^*(E)$-modules. If $n$ is even we have cofibre sequence
$$ C \lla C^*(B) \lla C^*(E)$$
and
$$ \cosusp{2n-2} C^*(B) \lla C \lla \cosusp{n-1}C^*(E)$$
of $C^*(E)$-modules, and by the octahedral axiom, the fibre of the composite
$$C^*(B) \lra C \lra \cosusp{2n-2} C^*(B)$$
is a $C^*(E)$-complex with one cell in codegree 0 and one cell in codegree
$n-1$. In particular, it is small.


\subsection{The situation in homology.}


We start with the case that $n$ is odd, since it  is a little simpler.

\begin{prop}
\label{prop:nodd}
If $n$ is odd, there is a short exact sequence
$$0 \lra \cosusp{n} H^*(B) \lra H^*(E)\lra H^*(B)\lra 0$$
\end{prop}

\begin{proof}
Since $p:E \lra B$ is split,  $H^*(E)$ is a free module over $H^*(B)$ on
two generators. We have a split surjection $s^*: H^*(E) \lra H^*(B)$,
so the kernel is a cyclic $H^*(B)$-module, and therefore a principal
ideal in $H^*(E)$, generated by an element $\tau$ of codegree $n$.
We therefore have a short exact sequence
$$0\lra (\tau) \lra H^*(E)\lra H^*(B)\lra 0, $$
and $(\tau ) \cong H^*(E)/\ann (\tau)$.
Since $n$ is odd, $\tau^2=0$, and hence $(\tau ) \subseteq \ann (\tau)$.
Since $H^*(E)$ is a free module on $1$ and $\tau$ it follows that
we have equality: $(\tau ) =\ann (\tau)$, giving the isomorphism
$$(\tau) \cong H^*(E)/\ann (\tau) \cong H^*(E)/(\tau)\cong H^*(B)$$
of $H^*(E)$-modules as required.
\end{proof}

If $n$ is even, it is easy to find examples where the conclusion
fails.
\begin{example}
Take $B=S^4$, $E=S^4\times S^4$, with $p$ being the projection onto
the first factor and $s$ the diagonal. Then
$H^*(E)=\Lambda (x_1,x_2)$, $H^*(B)=H^*(E)/(x_1-x_2)$ and

$$\ker (s^*)=\Sigma_n H^*(E)/(x_1+x_2)$$
so we do not have the conclusion of Proposition \ref{prop:nodd}
unless the characteristic is 2.
\end{example}

\begin{prop}
\label{prop:neven}
If $n$ is even, there are short exact sequences
$$0 \lra \cosusp{n} H^*(B)' \lra H^*(E)\lra H^*(B)\lra 0$$
and
$$0 \lra \cosusp{n} H^*(B) \lra H^*(E)\lra H^*(B)'\lra 0,$$
where $H^*(B)'$ is an $H^*(E)$-module isomorphic to $H^*(B)$ as
an $H^*(B)$-module.
\end{prop}

\begin{proof}
The proof begins as for Proposition \ref{prop:nodd}, but we need
no longer have $\tau^2=0$. To control the situation we use
the fact that the map $s^*: H^*E \lra H^*B$ is the edge homomorphism of the
Serre spectral sequence of the fibration
$$\Omega S^n \lra B \stackrel{s}\lra E,$$
and that $H^*(\Omega S^n)$ is the tensor product of an exterior algebra
on the generator $\Phi$ of codegree $n-1$ and the divided power algebra
on $\Psi$ of codegree $2n-2$. We have $d_n(\Phi)=\tau$ and we define $\tau_1$ by
$d_n(\Psi)=\tau_1\Phi$. It will be convenient to write
$$d_n^i: E_n^{*,i(n-1)} \lra E_n^{*,(i-1)(n-1)} $$
for the $H^*(E)$-module map giving the part of the differential
from the $i$th nonzero row to the $(i-1)$st.

By comparison with the path-loop fibration for $S^n$ we see that $\tau$
restricts to a generator of $H^n(S^n)$, and therefore, from the
Serre spectral sequence of $S^n \lra E \lra B$, that $1$ and $\tau$
give an $H^*(B)$-basis for $H^*(E)$.  This means that
$E_{n+1}^{*,0}=H^*(E) /(\tau) \cong H^*(B)$, and the edge homomorphism
shows that we have an exact sequence
$$0 \lra (\tau) \lra H^*(E) \stackrel{s^*}\lra H^*(B) \lra 0. $$
Since $H^*(E)$ is free as an $H^*(B)$-module on two generators,
we conclude $(\tau)$ is free as an
$H^*(B)$-module on one generator
of degree $n$, and we obtain the first exact sequence.

In particular $ (\tau)$ has a single copy of $k$ in the bottom
codegree.
Next, we have the exact sequence
$$0 \lra \ann (\tau) \lra H^*(E) \lra (\tau) \lra 0,  $$
and, since $H^*(E)$ is $H^*(B)$-free on two generators,
 $\ann (\tau)$ is free as an
$H^*(B)$-module on one generator
of degree $n$, and in particular it has a single copy of $k$ in the bottom
codegree. Of course $\ker (d_n^1)=\ann (\tau)$, whilst
on the other hand, by comparison with the path-loop fibration, $d_n^2(\Psi)$
generates a copy of $k$ in degree $n$, and since it lies in $\ker (d_n^1)$
we see $\tau_1$ generates $\ker (d_n^1)$, and
$$(\tau_1)=\ann(\tau).$$

Finally, we have the exact sequence
$$0 \lra \ann (\tau_1) \lra H^*(E) \lra (\tau_1) \lra 0. $$
Since $H^*(E)$ is $H^*(B)$-free on two generators,
 $\ann (\tau_1)$ is free as an $H^*(B)$-module on one generator
of degree $n$, and in particular it has a single copy of $k$ in the bottom
codegree. Of course $\ker (d_n^2)=\ann (\tau_1)$, whilst
on the other hand
$d_n(\Phi\Psi)=\tau \Psi + \Phi^2 \tau_1=\tau \Psi$ lies in this kernel and generates
generates a copy of $k$ in degree $n$. Since $d_n(\Phi\Psi)$  lies in $\ker (d_n^2)$,
we see $\tau$ generates $\ker (d_n^2)$ so that
$$(\tau)=\ann(\tau_1). $$
This identifies $(\tau_1)$ with $H^*(B)=H^*(E)/(\tau)$ as an $H^*(E)$-module,
giving the second exact sequence as required.
\end{proof}

\subsection{Lifting to cochains.}
\label{subsec:liftingtocochains}

Whether $n$ is even or odd, we have constructed a $C^*(E)$-module
$M$ for which $\pi_*(M)=H^*(B)$ as $H^*(E)$-modules. We now show
that we may lift this conclusion to the cochain level.

\begin{prop}
There is a unique $C^*(E)$-module $M$ with $\pi_*(M)\cong H^*(B)$ as
$H^*(E)$-modules.
\end{prop}

\begin{proof}
We will first give the proof assuming there is a short exact sequence
$$0 \lra \cosusp{n} H^*(B) \lra H^*(E)\lra H^*(B)\lra 0 $$
of $H^*(E)$-modules.

The first step is to show that $M$ is equivalent to a module
constructed by adding free $C^*(E)$-cells in codegrees
$0,n-1, 2n-2, 3n-3, \ldots$ as one might hope.
We take $M_0=M$, and then construct a diagram
$$M =M_0\lra M_1\lra M_2\lra \cdots $$
where
$$\pi_*(M_i)=\cosusp{(n-1)i} H^*(B).$$
To construct $M_{i+1}$ from $M_i$ we construct a cofibre sequence
$$\cosusp{ni-i} C^*(E) \stackrel{a_i}\lra M_i\lra M_{i+1}$$
where $a_i$ is a chosen generator of $\pi_*(M_i)$. It follows from the
hypothesis that $\pi_*(M_{i+1})=\cosusp{(n-1)(i+1)}H^*(B)$ as required to
proceed.
Since each map $M_i \lra M_{i+1}$ is zero in homotopy,
$M_{\infty}=\hcl{i}M_i$ has zero homotopy and is thus contractible.

Now construct the dual tower by the cofibre sequences
$$
\begin{array}{ccccc}
\downarrow &&\downarrow&&\downarrow\\
M^i& \lra &M &\lra& M_i \\
\downarrow &&\downarrow&&\downarrow\\
M^{i+1}& \lra &M &\lra& M_{i+1} \\
\downarrow &&\downarrow&&\downarrow
\end{array}$$
Passing to direct limits we see
$$M^{\infty} := \hcl{s} M^s \stackrel{\simeq}\lra M$$
is an equivalence. On the other hand,
 we have $M^0\simeq 0$ and by the octahedral axiom, we have
cofibre sequences
$$M^{i}\lra M^{i+1}\lra \cosusp{(n-1)i}C^*(E)$$
the module $M^{\infty}$ is a version of $M$ built with periodic cells $C^*(E)$
as required.

To see that any two such modules $M(0)$ and $M(1)$ are equivalent, we
perform the above construction to replace the two modules $M(a)$ ($a=0,1$)
with  filtered modules $\{ M(a)^i\}_i$
with $M(a)^{i+1}/M(a)^i=\cosusp{(n-1)i}C^*(E)$.
We suppose given an isomorphism $f_*: \pi_*(M(0))\lra \pi_*(M(1))$, and
we show it can be realized with a map $f:M(0) \lra M(1)$.
Indeed, we recursively construct maps
$$f_i: M(0)^i\lra M(1)$$
agreeing with the map $M(0)^i \lra M(0)$ in homotopy, where we identify
$\pi_*(M(0))$ and $\pi_*(M(1))$ using $f_*$.
Since $M(0)^0=0$, there is nothing to prove for $i=0$. After that
we are faced with the extension problem
\[ \xymatrix{
\cosusp{(n-1)i-1}C^*(E)\dto^{b_i}\\
M(0)^i\dto \rto^{f_i} &M(1)\\
M(0)^{i+1}\urdotted|>\tip_{f_{i+1}}&
} \]
Since the map $\cosusp{(n-1)i-1}C^*(E)\stackrel{b_i}\lra M(0)^i\lra M(0)$
is zero in
homotopy by construction, the same is true for the map
$\cosusp{(n-1)i-1}C^*(E)\stackrel{b_i}\lra M(0)^i\stackrel{f_i}\lra M(1). $
Since the domain is free, this shows $f_ib_i\simeq 0$,
and we can solve the extension problem. Passing to limits,
the Milnor exact sequence gives a map $M(0) \lra M(1)$ inducing a
homotopy isomorphism. This completes the proof if $n$ is odd.

If $n$ is even, the proof is precisely similar, except that the odd image
modules are suspensions of $H^*(B)'$ rather than of the standard
$H^*(E)$-module $H^*(B)$.
\end{proof}

\newcommand{\End}{\mathrm{End}}     

\section{Growth conditions}
\label{sec:eciisgci}
In this section we prove perhaps the simplest implication between the ci
conditions: for spaces of finite type, eci impies gci.

\subsection{Polynomial growth.}
Throughout algebra and topology it is common to use the rate of growth of
homology groups as a measurement of complexity. We will be working over
$H^*(X)$, so it is natural to assume that our modules $M$
are  {\em locally finite} in the sense that  $H^*(M)$ is cohomologically
 bounded below  and
 $\dim_k (H^i(M))$ is finite for all $i$.

\begin{defn}
We say that a locally finite module $M$ has polynomial growth of
degree $\leq d$ and write
$\growth (M) \leq d$ if
there is a polynomial $p(x)$ of degree $d$ with
$$\dim_k (H^{n}(M))\leq p(n)$$
for all $n >> 0$.
\end{defn}

\begin{remark}
Note that a complex with bounded homology has growth $\leq -1$. For
complexes with growth $\leq d$ with $d\geq 0$, by adding a constant to the
polynomial, we may insist that the
bound applies for all $n \geq 0$.
\end{remark}

\subsection{Mapping cones reduce degree by one.}
We give the following estimate on growth.
\begin{lemma}
\label{lem:reducegrowth}
Given locally finite modules $M,N$ in a triangle
$$\Sigma_n M \stackrel{\chi}\lra M \lra N$$
with $n \neq 0$ then
$$\growth (M) \leq \growth (N)+1.$$
\end{lemma}

\begin{proof}
The homology long exact sequence of the triangle includes
$$\cdots \lra H^{i-n}(M) \stackrel{\chi}\lra
H^i(M) \lra H^i(N)\lra \cdots .$$
This shows
$$\dim_k (H^i(M))\leq \dim_k (H^i(N))+\dim_k (\chi H^{i-n}(M)) .$$
Iterating $s$ times, we find
\begin{multline*}
\dim_k (H^i(M))\leq \dim_k (H^i(N))+\dim_k(H^{i-n}(N)) +\cdots
\\ \cdots +
\dim_k(H^{i-(s-1)n}(N)) +
\dim_k(\chi^s H^{i-sn}(M)) .
\end{multline*}

To obtain growth estimates, it is convenient to collect
the dimensions of the homogeneous parts into the Hilbert
series $h_M(t)=\sum_n\dim_k (H^i(M))t^i$. An inequality
between such formal series means that it holds between all
coefficients.

First suppose that $n>0$. Since  $H^*(M)$ is bounded below,
if $h_M(t)$ is the Hilbert series of $H^*(M)$ then we have
$$h_M(t) \leq h_N(t)(1+t^{n}+t^{2n}+\cdots )=\frac{h_N(t)}{1-t^{n}},$$
giving the required growth estimate.

If $n=-n'<0$ we rearrange to obtain
$$N'\lra M \lra \Sigma_{n'}M$$
where $N'=\Sigma_{n'-1}N$ and argue precisely similarly.
\end{proof}

\subsection{Growth of eci spaces.}
The implication we require is now straightforward.

\begin{thm}
If $X$ is eci  then it is also gci, and if $X$ has e-codimension $c$
it has g-codimension $\leq c$.
\end{thm}

\begin{proof}
It is sufficient to show $C^*(\Omega X)\simeq k \otimes_{C^*(X)}k $ has
polynomial growth.

By hypothesis there are elements $\chi_1,\chi_2, \ldots , \chi_c$
of non-zero degree in  $ZD(C^*(X))$ so that
$k/\chi_c/\ldots /\chi_1$ is small, and thus, applying
$k \otimes_{C^*(X)}(\cdot)$ we obtain a
complex with  growth $\leq -1$.
By Lemma \ref{lem:reducegrowth} if we apply  $k \otimes_{C^*(X)}(\cdot)$ to
$k/\chi_c /\ldots /\chi_2$ we obtain a complex
of growth $\leq 0$. Doing this repeatedly, we  deduce that when we
apply $k \otimes_{C^*(X)}(\cdot)$ to $k$ itself we obtain a complex
with  growth $\leq c-1$ as required.
\end{proof}

\section{Properties of gci spaces}
\label{sec:gciisGor}

The main results here are a smallness condition for modules over the
cochains of a normalizable gci space $X$ and a structural result for
the homology of $\Omega X$ which generalizes results of F{\'e}lix,
Halperin and Thomas. Throughout this section we assume the LRC \ref{LRC},
alternatively the reader may forego assuming the LRC and instead replace
the normalizability condition with polynomial normalizability.



\subsection{Cellular approximation, completion and Gorenstein condition}
We have need of recalling concepts from~\cite{DGI1} and~\cite{tec}.

Let $R=C^*(X)$. A map $U \to V$ of $R$-modules is a \emph{$k$-equivalence}
if the induced map $\Hom_R(k,U) \to \Hom_R(k,V)$ is an equivalence.
A module $N$ is \emph{$k$-null} if $\Hom_R(k,N)\simeq 0$. A module $C$ is
\emph{$k$-cellular} if $\Hom_R(C,N)\simeq 0$ for every $k$-null module $N$.
A module $C$ is \emph{$k$-complete} if $\Hom_R(N,C)\simeq 0$ for every
$k$-null module $N$. A map $f:C \to X$ of $R$-modules is a
\emph{$k$-cellular approximation (of $X$)} if $C$ is $k$-cellular and
$f$ is a $k$-equivalence. Finally, a map $f:X \to C$ is a
\emph{$k$-completion (of $X$)} if $C$ is $k$-complete and
$f$ is a $k$-equivalence.

We shall not recall the general definition of the Gorenstein condition
from~\cite{DGI1}; instead we give a definition which is equivalent under
the mild assumption of proxy-smallness.

Recall that an object is {\em proxy-small} if the thick category it
generates includes a non-trivial small object building the original object. Suppose that $k$ is
proxy-small as an $R$-module. Then $R$ is \emph{Gorenstein} if
$\Hom_R(k,R)\sim \Sigma^a k$ as left $k$-modules for some $a$ \cite[Prop. 8.4]{DGI1}.

In the situation we consider, where $X$ is normalizable, then $k$ is
always proxy-small by the following lemma.

\begin{lemma}
Let $X$ be a connected space. If $H^*(X)$ is
Noetherian then $k$ is proxy-small over $C^*(X)$.
\end{lemma}
\begin{proof}
Since $H^*(X)$ is Noetherian, it is finitely generated as a module
over a polynomial subring. We may use the polynomial generators to
form a Koszul complex $K$ on $C^*(X)$. Since its homology is finite
dimensional over $k$, it is finitely built by $k$. Finally, it is easy
to see that $K \otimes_{C^*(X)} k$ is non-zero. Since $k$ is a retract
of $K \otimes_{C^*(X)} k$, we see that $K$ builds $k$.
\end{proof}

\subsection{The Gorenstein property of gci spaces}
It is well known that for commutative rings a complete
intersection is Gorenstein. We prove an analogue, using the weakest of
the analogues of ci. 

\begin{prop}
\label{prop:gciisGorenstein}
If $X$ is a normalizable gci space then both $\chains_*(\Omega X)$ and
$\chains^*(X)$ are Gorenstein.
\end{prop}

In effect we rely on \cite{FHT} for the simply connected finite case,
and then apply the fibration lemma for Gorenstein spaces
\cite[10.2]{DGI1}  twice to deduce it in general. We
will present the proof in stages. Recall that $k=\Fp$.

\begin{lemma}
Let $F$ be a $k$-finite simply connected gci space,
then $\chains^*(F)$ is Gorenstein.
\end{lemma}

\begin{proof}
We begin by observing that $F$ is equivalent to the $p$-completion of a finite simply
connected CW-complex $F'$, which is thus a space of finite LS
category. Since $C_*(\Omega F')\simeq C_*(\Omega F)$, $F'$ is also
gci and therefore we see that $H_*(\Omega F')$ is an elliptic Hopf
algebra in the sense of \cite{FHT}. By \cite[Prop. 3.1]{FHT}, it is
Gorenstein in the sense that $\ext_{H_*(\Omega F')}^*(k,H_*(\Omega
F'))$ is one dimensional (non-zero in degree $a$ say), and the same is
true with $F'$ replaced by $F$.
Accordingly $\Hom_{C_*(\Omega F)}(k,C_*(\Omega F))\simeq \Sigma^a
k$ and $C_*(\Omega F)$ is Gorenstein in the sense of \cite{DGI1}.

By Proposition \ref{prop:EM}, the Eilenberg-Moore spectral sequence
for $C^*(F)$ converges,  so it follows from
\cite[Prop. 8.5]{DGI1} that $C^*(F)$ is also Gorenstein,
as required.
\end{proof}

Next we deal with the case of a finite complex which need not be
simply connected.

\begin{lemma}
\label{lem:FGor}
Let $F$ be a connected $k$-finite gci space with finite fundamental group such that
its universal cover $\tilde{F}$ is Gorenstein. Then $\chains^*(F)$ is Gorenstein.
\end{lemma}
\begin{proof}
Consider the fibration $\tilde{F} \lra F \lra B\pi_1F$. Since
$\pi_1(F)$ is a $p$-group, $\chains^*(B\pi_1(F))$ is dc-complete in
the sense of \cite[4.16]{DGI1} and so   $\chains^*(B\pi_1(F))$ is Gorenstein - because $\chains_*(\Omega B \pi_1(F))$ is. Since $\tilde{F}$ has finite homology, $\chains^*(F)$ is a small $\chains^*(B\pi_1(F))$-module. Now we can apply \cite[Prop. 8.10]{DGI1}, which shows that $\chains^*(F)$ is Gorenstein.
\end{proof}

Finally we are ready to deal with the general case.

\begin{proof}[ of Proposition \ref{prop:gciisGorenstein}]
First, suppose $F \lra X \lra B$ is a polynomial normalization for $X$.
By \cite[10.2]{DGI1} $\chains^*(B)$ is Gorenstein. 

\begin{lemma}
\label{fibreisgci}
The Noether fibre $F$ is gci.
\end{lemma} 

\begin{proof}
 First note that  both $F$ and $\Omega B$ are connected since the normalization 
is polynomial. Since $\Omega B$ is a finite $H$-space it is gci.

Next, we may assume that  $X$ is simply connected. Indeed,  
we may replace $X$ by its universal cover and $F$ by
its corresponding cover and  since $\pi_1(X)$ is finite, this will not 
affect the conclusion.  Now consider the fibration
$\Omega^2 B \lra \Omega F \lra \Omega X$. It  is principal
and the $\pi_1(\Omega X)$ action on $H_*(\Omega^2 B)$ is trivial.
The fact that $\Omega F$ has polynomial growth follows from the 
Serre spectral sequence. 
\end{proof}

It now follows from Lemma \ref{lem:FGor} that the Noether fibre $F$ is Gorenstein. Since $F$ is $k$-finite,
$\chains^*(X)$ is small as a $\chains^*(B)$-module. As before,
we apply \cite[Prop. 8.10]{DGI1}, which shows that
$\chains^*(X)$ is Gorenstein. Finally, since $\chains^*(X)$ is
dc-complete \cite[4.22]{DGI1},
then $\chains_*(\Omega X)$ is also Gorenstein, by \cite[Prop. 8.5]{DGI1}.

This completes the proof of Proposition \ref{prop:gciisGorenstein}
\end{proof}

\subsection{Completion and smallness criteria}
We shall need criteria for both smallness and $k$-completeness of
$C^*(X)$-modules.

\begin{lemma}
\label{lem: Completion over cochains of gci space}
Suppose $X$ is a normalizable gci space.
Then the $k$-completion of any $\chains^*(X)$-module $M$ is given by
\[ \Hom_{\chains_*(\Omega X)}(\Sigma^a k, \Hom_{\chains^*(X)}(k,M)) \]
where $a$ is the Gorenstein shift of $\chains^*(X)$.
\end{lemma}
\begin{proof}
Set $R=\chains^*(X)$ and $\ee=\chains_*(\Omega X)$.
Since $R$ is dc-complete, $\ee \simeq \Hom_R(k,k)$.
As we saw above, $k$ is proxy-small as an $R$-module,
so by \cite{DGI1}
\[ \cell_k R \simeq \Hom_R(k,R) \otimes_\ee k\]
Using the proxy-smallness of $k$ once again, it is
an easy exercise to show that the $k$-completion
of any $R$-module $M$ is given by
\[ M^\wedge_k \simeq \Hom_R(\cell_k R, M)\]
Bearing in mind that $k$ is an $R \otimes \ee$-module,
we have the following equivalences given by standard adjunctions:
\begin{align*}
\Hom_R(\cell_k R, M) & \simeq \Hom_R(\Hom_R(k,R) \otimes_\ee k, M) \\
& \simeq \Hom_{\ee^\op}(\Hom_R(k,R),\Hom_R(k,M))
\end{align*}

Since $X$ is Gorenstein,  we see that
$\Hom_R(k,R) \simeq \Sigma^a k$ as left $k$-modules.
Because $\pi_1(X)$ is a $p$-group and $k=\Fp$
there is only one simple $k[\pi_1(X)]$-module which is the trivial module $k$.
Hence the $\pi_0(\ee)$-module $\pi_a(\Hom_R(k,R))$ is that simple module $k$.
Now we can employ \cite[Prop. 3.9]{DGI1} which shows that
$\Hom_R(k,R) \simeq \Sigma^a k$ as right $\ee$-modules
\end{proof}

It is useful to have a criterion for $k$-completeness.
\begin{lemma}
If $X$ is a Gorenstein space  with $H^*(X)$ Noetherian then every
bounded-above $\chains^*(X)$-module of finite type
is $k$-complete.
\end{lemma}

\begin{remark}
(i) By Proposition~\ref{prop:gciisGorenstein} this applies to any
normalizable gci space.

(ii) The exact generality of the lemma is not clear, but we note that
$k$ need not be $k$-complete in general.  For example if $R=S^0$ and 
$k=H{\Bbb{F}}_p$ then $R$ is itself $k$-null by Lin's
theorem \cite{Lin}, and hence  $k$ is not $k$-complete.  
\end{remark}

\begin{proof}
Denote by $R$ the $k$-algebra $C^*(X)$. We first note that it is
sufficient to show that $k$ is $k$-complete. Indeed, in that case 
every $R$-module finitely built from $k$ is $k$-complete.
Since there are dual Postnikov sections in the category of $R$-modules
(see \cite[Prop. 3.3]{DGI1}) then every bounded-above $R$-module of finite type
is the homotopy limit of $k$-complete modules and hence is itself $k$-complete.

Since $k$ is proxy-small, there is a small object $K$ so that 
$$k\finbuilds K \builds k .$$
Because $R$ is Gorenstein (Proposition~\ref{prop:gciisGorenstein}) and
coconnective with $\pi_0(R)$ a field we have that
$\Hom_R(k,R)\simeq \Sigma^a k$ as a left $R$-module
(see~\cite[Prop. 3.9]{DGI1}). 

Now suppose $N$ is $k$-null, so that $\Hom_R(k,N)\simeq 0$.
For an $R$-module $M$ we denote by $DM$ the $R$-module $\Hom_R(M,R)$. 
Since $k\builds K$ we find $Dk \builds DK$, so that by the 
Gorenstein condition $k \builds DK$. Hence
$$0\simeq \Hom_R(DK,N)\simeq DDK \tensor_R N\simeq K \tensor_RN.$$
Since $K\builds k$ it follows that $k\tensor_R N\simeq 0$. 
Finally, 
\begin{align*}
\Hom_{R}(N,k) &\simeq \Hom_{R}(N,\Sigma^n \Hom_{R}(k,R))\\
&\simeq \Hom_R(k\otimes_R N,\Sigma^n R) \simeq 0
\end{align*}
and $k$ is $k$-complete as required. 
\end{proof}



We get the following criteria for smallness.
\begin{cor}
\label{cor: Condition for smallness}
Suppose $X$ is a normalizable gci space and let $M$ be a bounded-above
$\chains^*(X)$-module of finite type. The following are equivalent
\begin{enumerate}
\item $M$ is small,
\item $\pi_*(M\otimes_{\chains^*(X)} k)$ is finite.
\end{enumerate}
\end{cor}
\begin{proof}
Set $R=\chains^*(X)$. Clearly if $M$ is small then
$\pi_*(M\otimes_R k)$ is finite.

So suppose that $M \otimes_R k$ has finite homotopy.
Let $\ee=\chains_*(\Omega X)$. Since $k$ is proxy-small
and $M$ is $k$-complete we see that $M=\Hom_\ee(\Sigma^a k,\bar{M})$
where $\bar{M}$ is the $\ee$-module $\Hom_R(k,M)$.

By \cite[Prop. 4.17]{DGI1} $k$ is proxy-small also as
an $\ee$-module. It follows that $\Sigma^a M\otimes_R k$ is a
$k$-cellular approximation of $\bar{M}$. Hence there is a $k$-equivalence
of $\ee$-module $\Sigma^a M\otimes_R k \to \bar{M}$. In particular
\[ M \simeq \Hom_\ee(k, M\otimes_R k)\]

Since $\pi_1(X)$ is a finite $p$-group it follows that
$\pi_n(M \otimes_R k)$ is finitely built by $k$ as a $k[\pi_1 X]$-module
for every $n$.
Using \cite[Prop. 3.2]{DGI1} and the fact that $M \otimes_R k$ has
finite homotopy we see that $k \finbuilds_\ee M\otimes_R k$ and therefore
\[ R\sim \Hom_\ee(k,k) \finbuilds_R \Hom_\ee(k,M\otimes_R k) \sim  M\]
\end{proof}

\subsection{Loop space homology of gci spaces}
The following rests heavily on results of F{\'e}lix, Halperin and Thomas.
\begin{lemma}
\label{lem: Loop homology of Gorenstein gci}
Let $X$ be a Gorenstein gci space. Then
\begin{enumerate}
\item $H_*(\Omega X)$ is left and right Noetherian and
\item $H_*(\Omega X)$ is a finitely generated module over a central polynomial sub-algebra.
\end{enumerate}
\end{lemma}
\begin{proof}
First suppose that $X$ is simply connected. The depth of the Hopf algebra $H_*(\Omega X)$ is the least integer $m$ such that $\ext^m_{H_*(\Omega X)}(k,H_*(\Omega X))$ is non-zero. We claim that $H_*(\Omega X)$ has finite depth. Observe that there is conditionally convergent spectral sequence
\[ E_2^{p,q}=\ext^{p,q}_{H_*(\Omega X)}(k,H_*(\Omega X)) \ \Rightarrow \ \ext^{p+q}_{\chains_*(\Omega X)}(k,\chains_*(\Omega X)) \]
Since $\chains_*(\Omega X)$ is Gorenstein, the $E_2$ term of this spectral sequence cannot be zero. Hence $H_*(\Omega X)$ has finite depth. Now, by Theorem C and Theorem B of \cite{FelixHalperinThomasHopfAlgebraElliptic}, $H_*(\Omega X)$ has the desired properties.

Now suppose that $X$ is not simply connected. Let $\tilde{X}$ be the
universal cover of $X$. Clearly $\tilde{X}$ is gci.
Using \cite[Prop. 8.10]{DGI1} it follows
from the fibration $\pi_1(X) \lra \tilde{X} \lra X$ that $\tilde{X}$ is
also Gorenstein. Hence $H_*(\Omega \tilde{X})$ has the desired
properties, i.e., it is left and right Noetherian and is a finite module over a central polynomial subalgebra.

There is a natural way of identifying $H_*(\Omega X)$ with the semi-direct product $H_*(\Omega \tilde{X}) \rtimes \pi_1(X)$. Since $\pi_1(X)$ is a finite group, we see that $H_*(\Omega X)$ is a left and right Noetherian  $H_*(\Omega \tilde{X})$-module. In particular $H_*(\Omega X)$ is left and right Noetherian (over itself).

Note that $\pi_1(X)$ acts on $H_*(\Omega\tilde{X})$ via algebra maps and
therefore the center of $H_*(\Omega\tilde{X})$, which we shall denote by
$\tilde{Z}$, is invariant under this action. Because there is a
polynomial subalgebra $P \subset \tilde{Z}$ such that $H_*(\Omega
\tilde{X})$ is a finitely generated $P$-module, $\tilde{Z}$ is
Noetherian. The Hilbert-Noether theorem shows
that the ring $\tilde{Z}^{\pi_1(X)}$ of invariants is Noetherian and
$\tilde{Z}$ is a finitely generated over it. We conclude that
$H_*(\Omega X)$ is finitely generated over $\tilde{Z}^{\pi_1(X)}$.

From the identification
$H_*(\Omega X) \cong H_*(\tilde{X}) \rtimes \pi_1(X)$ we see that
$\tilde{Z}^{\pi_1(X)}$ is contained in the center of $H_*(\Omega X)$.
Therefore $H_*(\Omega X)$ is finitely generated over its center $Z$,
and $Z$ is Noetherian.
From the Noether normalization theorem we conclude that
$H_*(\Omega X)$ is a finitely generated module over a
central polynomial subalgebra.
\end{proof}

\section{The Hochschild cohomology shearing map}
\label{sec:shearing}

We need to discuss  certain ring homomorphisms
$$\psi : HH^*(X) \lra H_*(\Omega X)$$
in the absolute case and
$$\chi : HH^*(X|B) \lra H_*(\Omega F)$$
when we have a fibration $F\lra X \lra B$. We refer to these as
{\em shearing maps}, and  the purpose of the present section is to show that
several different possible definitions agree. Note that $\psi$ is well known
in various contexts, see e.g. \cite{FTVP}.

\subsection{General context}
In this section we give two equivalent definitions for the shearing
map. Throughout we assume $Q$ is a commutative ring spectrum and that
there are maps of ring spectra
\[ Q \xrightarrow{\rho} R \xrightarrow{\sigma} S\]
Thus, $R$ and $S$ are $Q$-algebras, note that we do not assume $R$ and $S$ are commutative $Q$-algebras. This is a precaution, as it is often all too easy to assume the wrong bimodule structure when working in a commutative setting.

In this section we shall denote by $R^e_Q$ the $Q$-algebra
$R\otimes_Q R^\op$. We will denote by $S^\sigma$ the $Q$-algebra $R\otimes_Q S^\op$. Note there are maps of $Q$-algebras
\[ R^e_Q \xrightarrow{1\otimes \sigma} S^\sigma \xrightarrow{\sigma \otimes 1} S^e_Q\]
Since $S$ is an $S^e_Q$-module, this structure is pulled back to make $S$ both an $S^\sigma$-module and an $R^e_Q$-module. The shearing map is a map of graded algebras
\[ HH^*(R|Q) \longrightarrow \ext^*_{S^\sigma} (S,S)\]

This is probably the place to mention we have two main settings in mind. In the first there is a fibration $F \to X \to B$ and we set $Q=C^*(B)$, $R=C^*(X)$ and $S=C^*(F)$ with $\rho$ and $\sigma$ being the obvious maps induced by the fibration. In this first setting $S^\sigma$ is equivalent to the $k$-algebra $S\otimes_k S$ and so $\ext^*_{S^\sigma} (S,S)$ is the Hochschild cohomology $HH^*(F)$. In the second setting $R=C^*(X)$ for some space $X$ and both $S$ and $Q$ are $k$. Here $S^\sigma$ is equivalent to $R$ and $\ext^*_{S^\sigma}(S,S)$ turns out to be $H_*(\Omega X)$. There is also a third setting, which is an amalgamation of the first two. We shall say more on all of these settings towards the end of this section.

\subsection{Two descriptions of the shearing map}

Before describing the shearing map we need a pair of lemmas. In what follows we refer to the right copy of $R$ in $R^e_Q$ as $R_r$, thus there is a $k$-algebra map $R_r \to R^e_Q$.

\begin{lemma}
\label{lem: Functor to S beta modules}
There is a functor $F:\Derived(R^e_Q) \to \Derived(S^\sigma)$ given by $M \mapsto M\otimes_{R_r} S$. This functor is naturally equivalent to the functor $M \mapsto S^\sigma \otimes_{R^e_Q} M$ and therefore $F$ is left adjoint to the forgetful functor $\Derived(S^\sigma) \to \Derived(R^e_Q)$.
\end{lemma}


\begin{lemma}
\label{lem: Equivalence of S beta modules}
There is an equivalence of $S^\sigma$-modules: $R \otimes_{R_r} S \simeq S$.
\end{lemma}

At first glance it is perhaps not clear what is the claim in Lemma~\ref{lem: Equivalence of S beta modules} above. Recall that $S$ is an $S^e_Q$-module, and it is this module structure that we pulled back to make $S$ into an $S^\sigma$-module. So it is not immediately apparent that the $S^\sigma$-module structure on $R \otimes_{R_r} S$ agrees with the one on $S$.

The proofs of Lemmas~\ref{lem: Functor to S beta modules} and~\ref{lem: Equivalence of S beta modules} are both based on chasing the same diagram. There is a commuting square of maps of $Q$-algebras
\[ \xymatrixcompile{
{R_r^\op} \ar[r] \ar[d] & {S^\op} \ar[d] \\
{R^e_Q} \ar[r] & {S^\sigma} }
\]
This induces a diagram of homotopy categories and adjoint functor between
them
\[ \xymatrixcompile{
{\Derived(R_r^\op)} \ar@<1ex>[r] \ar@<1ex>[d] & {\Derived(S^\op)} \ar@<1ex>[l] \ar@<1ex>[d] \\
{\Derived(R^e_Q)} \ar@<1ex>[r] \ar@<1ex>[u] & {\Derived(S^\sigma)} \ar@<1ex>[l] \ar@<1ex>[u] }
\]
Chasing the diagram above will easily yield proofs for Lemmas~\ref{lem: Functor to S beta modules} and~\ref{lem: Equivalence of S beta modules}, and we leave it to the reader to complete the details.

\begin{defn}
\label{def: Shearing map}
Define the \emph{shearing map for the transitivity triple $Q \to R \to
  S$} to be the map $\shr:\ext^*_{R^e_Q}(R,R) \to
\ext^*_{S^\sigma}(S,S)$ of graded rings given by $\shr(f) = f \otimes_{R_r} S$.
\end{defn}

By Lemmas \ref{lem: Functor to S beta modules} and \ref{lem: Equivalence of S beta modules} there is an isomorphism $\tau_R:\ext^*_{R^e_Q}(R,S) \to \ext^*_{S^\sigma}(S,S)$. This allows us to construct another map
\[ \alpha:\ext^*_{R^e_Q}(R,R) \to \ext^*_{S^\sigma}(S,S) \]
where $\alpha$ is the composition
\[ \ext^*_{R^e_Q}(R,R) \xrightarrow{\sigma} \ext^*_{R^e_Q}(R,S) \xrightarrow{\tau_R} \ext^*_{S^\sigma}(S,S)\]
We will show that $\alpha$ and $\shr$ are equal, but for that we must first have an explicit description of $\tau$.

Since $\tau$ is the isomorphism given by an adjunction, it only
depends on a choice of counit map $m:S\otimes_{R_r} S \to S$. This
then becomes a bootstrapping problem, namely describing the counit
map. For that we state the following lemma.

\begin{lemma}
\label{lem: The counit map}
There is the following commutative diagram of $S^e_Q$-modules
\[ \xymatrixcompile{ {S \otimes_Q S} \ar[r]^-{m'} \ar[d] & S \ar[d]^=\\ {S\otimes_R S} \ar[r]^-m & S} \]
where $m'$ is the multiplication map of $S$ as a $Q$-algebra. The map $m$ satisfies the following identities:
\begin{align*}
m(\sigma \otimes_R S) &= 1_S \\
m(S \otimes_R \sigma) &= 1_S
\end{align*}
where the first is an identity of morphisms in $\Derived(S^\sigma)$ and the second is an identity of morphisms in $\Derived(S \otimes_Q R^\op)$.\qqed
\end{lemma}

\begin{remark}
\label{rem: concerning the multiplication of S over R}
Suppose $F \lra X$ is a fibration and set $Q=k$, $R=C^*(X)$ and $S=C^*(F)$. Then the map $m'$ is induced by the diagonal $\Delta':F \lra F \times F$ while $m$ is induced by the diagonal $\Delta:F \lra F\times_X F$. Clearly $F\times_X F \subseteq F\times F$ and the composition $F \stacklra{\Delta} F\times_X F \lra F\times F$ is $\Delta'$. This gives a topological explanation to the commuting diagram above.
\end{remark}

Denote by $F$ the functor $-\otimes_{R_r} S:\Derived(R^e_Q) \to \Derived(S^\sigma)$ and by $G$ the forgetful functor $\Derived(S^\sigma) \to \Derived(R^e_Q)$. Let $\epsilon:FG \to 1$ be given by $\epsilon_Y=Y\otimes_S m$ and let $\eta:1 \to GF$ be $\eta_X=X\otimes_{R_r} \sigma$. From Lemma~\ref{lem: The counit map} we get the following identities:
\begin{align*}
1_{FX}&=\epsilon_{FX} \circ F(\eta X) \\
1_{GY}&=G(\epsilon Y) \circ \eta_{GY}
\end{align*}
Accordingly,  $\epsilon$ and $\eta$ are the counit and the unit for the adjunction of $F$ and $G$. In particular we get the following corollary.

\begin{cor}
\label{cor: Natural transformations for shearing map}
Define natural transformations $\tau:\ext^*_{R^e_Q}(-,S) \to \ext^*_{S^\sigma}(-\otimes_{R_r} S,S)$ and $\delta:\ext^*_{S^\sigma}(-\otimes_{R_r} S,S) \to \ext^*_{R^e_Q}(-,S)$ in the following manner. Given a map $f:M \to \Sigma^n S$ of $R^e_Q$-modules, let $\tau(f)$ be the composition $m(f \otimes_{R_r} S)$, and given a map $g:M\otimes_R S \to \Sigma^n S$ of $S^\sigma$-modules, let $\delta(g)=g(1_M \otimes_{R_r} \sigma)$. Then $\delta \tau = 1$ and $\tau \delta = 1$.
\end{cor}
\begin{proof}
The adjoint to a morphism $f:M \to GS$ is the composition $FM \xrightarrow{Ff} FGS \xrightarrow{\epsilon_S} S$. The adjoint to a morphism  $g:FM \to S$ is the composition $M \xrightarrow{\eta_M} GFM \xrightarrow{Gg} GS$. Since $\epsilon_S = m$ and $\eta_M=1 \otimes_{R_r} \sigma$, we are done.
\end{proof}


\begin{lemma}
\label{lem: nu and alpha are equal}
The map $\alpha$ defined above is equal to the shearing map $\shr$. In particular $\alpha$ is a map of graded algebras. Moreover, the right $\ext^*_{R^e_Q}(R,R)$-module structure on $\ext^*_{R^e_Q}(R,S)$ is the same as the module structure induced by the map of algebras $\shr$.
\end{lemma}
\begin{proof}
Let $f$ be an element of $\ext^*_{R^e_Q}(R,R)$, then
\begin{align*}
\tau\sigma (f)&=\tau(\sigma f)\\
&=m(\sigma f\otimes S)\\
&=m(\sigma \otimes S)(f\otimes S)\\
&=m(\sigma \otimes S)\shr(f)\\
&=\shr(f)
\end{align*}
The last identity comes from $m(\sigma \otimes S)=1_S$.

The second statement requires us to show that for $f\in \ext^*_{R^e_Q}(R,R)$ and $g \in \ext^*_{R^e_Q}(S,S)$:
\[ m(g\sigma f\otimes S)=g(f\otimes S)\]
Thus, it is enough to show that $m(g\sigma \otimes S)=g$, but this is precisely the statement $\tau \delta (g)=g$ which is proved above.
\end{proof}

\subsection{The shearing map for a fibration}
Let $F \stacklra{q} X \stacklra{p} B$ be a fibration and set $Q=C^*(B)$, $R=C^*(X)$ and $S=C^*(F)$. We shall show that $S^\sigma$ is equivalent to the $k$-algebra $S\otimes_k S$ and so $\ext^*_{S^\sigma} (S,S)$ is the Hochschild cohomology $HH^*(F)$.

\begin{lemma}
There is an equivalence of $k$-algebras $S^\sigma \simeq S\otimes_k S^\op=S^e$.
\end{lemma}
\begin{proof}
Recall that $S \simeq R \otimes_Q k$. Thus,
\begin{align*}
R^e_Q\otimes_R S^\op &\simeq (R\otimes_Q R) \otimes_R S^\op \simeq R\otimes_Q S^\op\\
&\simeq R\otimes_Q (k \otimes_Q R) \simeq (R\otimes_Q k)\otimes_k (k\otimes_Q R)\\
&\simeq S\otimes_k S^\op
\end{align*}
\end{proof}

\begin{remark}
There is also topological argument for the equivalence above. It is easy to show that the homotopy pullback of the diagram $X\times_B X \stacklra{p_2} X \xleftarrow{q} F$ is $F \times F$, where $p_2$ is the projection onto the second coordinate.
\end{remark}

\begin{remark}
The equivalence $S^\sigma \simeq S^e$ allows us to translate the
$S^\sigma$-module $R\otimes_R S$ into an $S^e$-module, which we shall
provisionally denote by $T$. Clearly $T$ is equivalent to $S$ as
$k$-modules, but one might worry that $T$ would not have the correct
$S^e$-module structure we want. The module structure we need is the
diagonal one, induced by the diagonal map $\Delta: F \lra F \times F$
(whereas, for example, the $S \otimes_k S$-module $k\otimes_k S$ has a
``wrong" module structure). We sketch a topological argument as to why $T$ has the correct module structure. Consider the map of homotopy pullback squares below, where the left and right walls are homotopy pullbacks and the map $X \times_B X \lra X$ is projection to the right coordinate:
\[ \xymatrixcompile{
{F} \ar[rr]^{\theta} \ar[dd] \ar[dr] & &
 {F \times F} \ar'[d]^(.7){}[dd] \ar[dr]^{} \\
& {X} \ar[rr]^\Delta \ar[dd] &&
 {X\times_B X} \ar[dd] \\
{F} \ar[dr] \ar'[r][rr] && {F} \ar[dr] \\
&{X} \ar[rr]^{=} &&{X}
} \]
One easily sees that the induced morphism $\theta$ of homotopy pullbacks is indeed the diagonal map, as required. Hence $T$ has the correct module structure.
\end{remark}

Thus we have shown.
\begin{cor}
In the setting above the shearing map
\[ \chi: HH^*(X|B) \lra HH^*(F)\]
 is a map of graded algebras.
\end{cor}

\subsection{The shearing map to the loop space homology}
\label{sub: The shearing map to the loop space homology}
In this setting $R=C^*(X)$ for some space $X$ and both $S$ and $Q$ are $k$. In this case $S^\sigma=R$ and $\ext^*_{S^\sigma}(S,S)\cong H_*(\Omega X)$. One immediately sees that the shearing map has the form
\[ \psi:HH^*(X) \lra H_*(\Omega X)\]
Note that~\cite{FTVP} describes a map $I:HH^*(X) \lra H_*(\Omega X)$ which is the same as the map $\alpha$ introduced earlier and so, by Lemma~\ref{lem: nu and alpha are equal}, it is the shearing map. It is shown in~\cite{FTVP} that the image of $I$ is central in $H_*(\Omega X)$, we give here a different argument for that fact.

In the current setting the shearing map is clearly equal to the composition
\[ HH^*(R) \lra ZD(R) \lra \ext^*_R(k,k)\]
where the map $ZD(R) \lra \ext^*_R(k,k)$ is given by $\zeta \mapsto \zeta_k:k \lra \Sigma^n k$. The following is now obvious.
\begin{lemma}
\label{lem:shearing image is central}
In this setting the image of the shearing map $\psi$ is contained in the graded centre of $H_*(\Omega X)$.
\end{lemma}


\subsection{The relative shearing map to the loop space homology}
\label{sub: Relative shearing map}
The setting here is as follows. Let $F \lra X \lra B$ be a fibration and set $Q=C^*(B)$, $R=C^*(X)$ and $S=k$. Then $S^\sigma=R\otimes_Q S$ turns out to be $C^*(F)$. Thus the shearing map is:
\[ \shr:HH^*(X|B) \lra H_*(\Omega F)\]

\begin{lemma}
\label{lem: composing with relative shearing map}
The composition
\[HH^*(X|B) \stacklra{\shr} H_*(\Omega F) \lra H_*(\Omega X)\]
is equal to the composition
\[HH^*(X|B) \lra HH^*(X) \stacklra{\psi} H_*(\Omega X)\]
where $\psi$ is the shearing map from \ref{sub: The shearing map to the loop space homology} above.
\end{lemma}
\begin{proof}
Let $a:R \lra k$ be the augmentation map and let $\varphi:R^e \lra R^e_Q$ be the obvious map. Note that $\varphi$ induces the forgetful functor $\varphi^*:\Derived(R^e_Q) \lra \Derived(R^e)$. In light of Lemma~\ref{lem: nu and alpha are equal} we need to show that the two compositions
\begin{eqnarray}
\ext^*_{R^e_Q}(R,R) \stacklra{a} \ext^*_{R^e_Q}(R,k) \stacklra{\varphi^*} \ext^*_{R^e}(R,k) &  \\
\ext^*_{R^e_Q}(R,R) \stacklra{\varphi^*} \ext^*_{R^e}(R,R) \stacklra{a} \ext^*_{R^e}(R,k) &
\end{eqnarray}
are equal.

Let $f:R \lra \Sigma^n R$ be an element of $\ext^*_{R^e_Q}(R,R)$. The first composition yields the morphism $\varphi^*(af) \in \ext^*_{R^e}(R,k)$ while the second composition yields $\varphi^*(a)\varphi^*(f) \in \ext^*_{R^e}(R,k)$. Both compositions are equal, because $\varphi^*$ is a functor.
\end{proof}

\begin{cor}
\label{cor: relative shearing map is composition of two others}
The shearing map $\shr:HH^*(X|B) \lra H_*(\Omega F)$ is equal to the composition
\[ HH^*(X|B) \stacklra{\alpha} HH^*(F) \stacklra{\beta} H_*(\Omega F)\]
where $\alpha$ and $\beta$ are the appropriate shearing maps.
\end{cor}
\begin{proof}
Given $f\in HH^*(X|B)$ we have
\[ \shr(f) = f\otimes_{C^*(X)} k = f \otimes_{C^*(X)}C^*(F) \otimes_{C^*(F)} k = \beta(\alpha(\shr))\]
\end{proof}

\section{Finite gci spaces are hhci}
\label{sec:gciiszci}

In this section we assume the LRC \ref{LRC}.

\subsection{A Sufficient condition for being hhci}
This condition (Theorem~\ref{thm: Sufficient condition} below)
will also be used in Section~\ref{sec:gciiseci}, and for that reason
we phrase it in greater generality than is
strictly needed in this section.

Throughout this section we denote by $\shr$ the shearing
map $HH^*(X|B) \lra H_*(\Omega F)$ where $X \lra B$ is a
polynomial normalization and $F$ is the Noether fiber.
Recall (Lemma~\ref{lem:shearing image is central})
that the image of $\shr$ is
contained in the graded-commutative centre of $H_*(\Omega F)$,
which we shall denote by $Z$. 
We also remind the reader that by Lemma~\ref{fibreisgci} the
space $F$ is gci whenever $X$ is gci.

\begin{lemma}
\label{lem: A is Noetherian}
Let $X$ be a normalizable gci space. Denote by $A$ the
image of the shearing
map $\shr$. Suppose that $Z$ is finitely generated as a module
over the subalgebra $A$. Then $A$ is a Noetherian algebra.
\end{lemma}
\begin{proof}
By Lemma~\ref{lem: Loop homology of Gorenstein gci} the algebra
$H_*(\Omega F)$ is a free finitely generated module over a central
polynomial subalgebra $P$. Let $Z\subseteq H_*(\Omega F)$ be the
graded-commutative centre. Then $A \subseteq Z$ and $Z$ is a
Noetherian graded-commutative ring ($Z$ is Noetherian because it is
also finitely generated as a $P$-module). In
\cite{EisenbudSubringsOfNoetherianRings} 
Eisenbud proves that if a central subring of a Noetherian ring
finitely generates the whole ring as a module, then that subring is
Noetherian. A simple generalization of this result to the
graded-commutative setting completes the proof. 
\end{proof}

\begin{thm}
\label{thm: Sufficient condition}
Let $X$ be a normalizable gci space and suppose that
the shearing map $\shr$ makes $Z$ into a finitely
generated $HH^*(X|B)$-module. Then $X$ is hhci.
\end{thm}
\begin{proof}
Let $A$ be the image of $\shr$. Since $A$ is a graded-commutative
Noetherian algebra (by Lemma~\ref{lem: A is Noetherian}) with
$A_0=k$, then $A$ has a Noether normalization. There are
elements $x_1,...,x_c \in A$ which generate a polynomial subalgebra
$P=k[x_1,...,x_c]$ such that $A$ is a finitely generated
$P$-module. Since $Z$ is finitely generated as an $A$-module,
it is also finitely generated as a $P$-module.

We know from Lemma~\ref{lem: Loop homology of Gorenstein gci} that
$H_*(\Omega F)$ is finitely generated as a $Z$-module and hence it
is also finitely generated as a $P$-module.

Set $R=C^*(X)$ and $Q=\chains^*(B)$ and choose elements
$z_1,...,z_c \in HH^*(R|Q)$ such that $\shr(z_i)=x_i$.
We will show that these elements satisfy the hhci definition.
Let $B=R/z_1/\cdots/z_c$ be the relevant $R^e_Q$-module;
we must show that $B$ is small.

Clearly $B$ is bounded-above and of finite type, thus we may
use the smallness criterion of
Corollary~\ref{cor: Condition for smallness}.
Note we are applying the criterion to modules over the
cochains of the space $X\times_B X$ - which is also a normalizable
gci space using the normalization $X\times_B X \lra B$.
Thus we must show that $B\otimes_{R^e_Q} k$ has finite homology.

Recall from Subsection~\ref{sub: Relative shearing map} that the transitivity
triple used for defining the shearing map $\shr$ is
$Q=\chains^*(B)$, $R=\chains^*(X)$ and $S=k$ and thus
$S^\sigma$ is $\chains^*(F)$.
As noted in Lemma~\ref{lem: Functor to S beta modules}
$B\otimes_R k \simeq S^\sigma \otimes_{R^e_Q} B$ and so
\[ B\otimes_{R^e_Q} k \simeq k \otimes_{S^\sigma} (B\otimes_R k).\]

Next we claim it is enough to show that
\[ \ext^*_{S^\sigma}(B\otimes_R k, k) \]
is finite dimensional. Indeed
\[ \Hom_{S^\sigma}(B\otimes_R k , k) \simeq
\Hom_k(k \otimes_{S^\sigma}(B\otimes_R k) , k)\]

Recall there is an isomorphism
$\ext^{*}_{S^\sigma}(k,k) \cong H_*(\Omega F)$.
For every $i=1,...,c$ the map
$z_i \otimes_R k \in \ext^*_{S^\sigma}(k,k)$
represents the element $x_i=\shr(z_i)$. From this it is
easy to see that the map
\[\ext^{*}_{S^\sigma}(z_i\otimes_R k,k):
H_*(\Omega F) \lra H_{*-|z_i|}(\Omega F)\]
is simply right multiplication by $x_i$.

The Koszul filtration on $B=R/z_1/\cdots/z_c$
induces a filtration on $B \otimes_R k$. This
gives rise to a spectral sequence whose $E^2$-term is the Koszul
homology
\[H_{s,t}(H_*(\Omega F)/x_1/\cdots/x_c) \cong
\tor^P_{s,t}(H_*(\Omega F), k)\]
which strongly converges to $\ext^{-s-t+c}_{S^\sigma}(k,B\otimes_R k)$.
Since $H_*(\Omega F)$ is a finitely generated $P$-module the $E^2$-term
of this spectral sequence is finite dimensional. Therefore
$\ext^*_{S^\sigma}(k,B\otimes_R k)$ is finite dimensional and
the proof is done.
\end{proof}

\subsection{Using a Hochschild cohomology spectral sequence}

The Hochschild cohomology spectral sequence of~\cite{ShamirHHSS2}
identifies one case where $H_*(\Omega F)$ is finitely generated
over the image of the shearing map.

\begin{thm}
\label{thm: Finite gci spaces are hhci}
A finite gci space is hhci.
\end{thm}
\begin{proof}
Let $X$ be a finite gci space. The polynomial normalization we
consider is of course $X \lra pt$ and the resulting shearing map
is $\shr:HH^*(X) \lra H_*(\Omega X)$.

Since $X$ is a normalizable gci space, then by
Lemma~\ref{lem: Loop homology of Gorenstein gci} the loop space
homology $H_*(\Omega X)$ is finitely generated over a central
polynomial sub-algebra. We can now use the results of~\cite{ShamirHHSS2}
which show that the centre of $H_*(\Omega X)$ is finitely generated over 
the image of the shearing map.
Invoking Theorem~\ref{thm: Sufficient condition} completes the proof.
\end{proof}

\section{Normalizable gci spaces are eci}
\label{sec:gciiseci}

Throughout $X$ is a normalizable gci space and we assume
the LRC \ref{LRC}. Our goal is to prove the following result.

\begin{thm}
\label{thm:gciiseci}
If $X$ is a normalizable gci space then $X$ is eci.
\end{thm}

The missing ingredient for showing that a normalizable gci space $X$
is hhci is described in the next result.

\begin{prop}
\label{pro: Normalizable gci is hhci sometimes}
Let $X$ be a normalizable gci space, let $\nu:X \lra B$ be a
polynomial normalization for $X$ and let $F$ be its Noether
fibre. If the relative Hochschild cohomology shearing map
$\shr:HH^*(X|B) \lra HH^*(F)$ makes $HH^*(F)$ into a finitely
generated $HH^*(X|B)$-module, then $X$ is hhci.
In particular, if $HH^*(X|B)$ is Noetherian then $X$ is hhci.
\end{prop}

\begin{example}
For two closely related examples for $X\lra B$ to which this applies
suppose $p$ is an odd prime. The two examples are 
(i) $BC_p\lra BSO(2)$ and (ii) $BD_{2p}\lra BO(2)$. In both cases
$H^*(B)=k[x]$ is polynomial on one generator (of degree 2 or 
4 respectively) and $H^*(X)=k[x]\tensor
\Lambda (\tau)$ where $|\tau|=|x|-1$. Since the spectral sequences
relating algebra and topology collapse  the shearing map is
$$HH^*(k[x]\tensor \Lambda (\tau)|k[x]) \lra HH^*(\Lambda (\tau)|k)$$
and  easily seen to be surjective. 
\end{example}

\subsection{Normalizable gci spaces}
Recall from Lemma~\ref{fibreisgci} that because $X$ is gci, 
then so is the Noether fibre $F$ of any polynomial 
normalization $X \lra B$. 
Theorem~\ref{thm: Finite gci spaces are hhci} above shows
that $F$ is hhci; we use this together with the normalizability
of $X$ to deduce it is eci.

\begin{proof}[ of Theorem~\ref{thm:gciiseci}]
Let $\nu:X \lra B$ be a polynomial normalization for $X$ and let $F$
be the fibre of this map, then $F$ is $k$-finite and is gci.

Set $Q=C^*(B)$, $R=C^*(X)$ and $S=C^*(F)$ and let
$x_1,...,x_n \in H^*(B)$ be the polynomial generators of
the cohomology algebra. So $k\simeq Q/x_1/\cdots/x_n$ and
\[ S \simeq R\otimes_Q k \simeq R/x_1/\cdots/x_m\]
as $R$-modules.

There is a commutative diagram of algebras
\[ \xymatrixcompile{
{R^e_Q} \ar[r] \ar[d] & {S^e} \ar[d] \\
{R} \ar[r]            & {S}
}\]
coming from the obvious diagram of spaces. Pulling back along the
left-hand vertical,  $R/x_1/\cdots/x_n \simeq S$ as $R^e_Q$-modules. Thus we have
shown that $R \finbuilds_{hh} S$.

By Theorem~\ref{thm: Finite gci spaces are hhci} the Noether fibre
$F$ is hhci. In particular there are elements
$z_1,...,z_c\in \ext^*_{S^e}(S,S)$ such that $T=S/z_1/\cdots/z_c$ is a
small $S^e$-module. Using the diagram above we can pull this
construction back to the category of $R^e_Q$-modules, thereby showing
that $S \finbuilds_e T$ as $R^e_Q$-modules. Combined with
the fact that $R \finbuilds_{hh} S$ this yields 
$R \finbuilds_e T$.

It remains to show that $T$ is small as an $R^e_Q$-module.
The homotopy fibre of the map $F \times F \lra X \times_B X$
is easily seen to be $\Omega B$. This implies that $S^e$
is small as an $R^e_Q$-module and therefore $T$ is small
as an $R^e_Q$-module.
\end{proof}

\subsection{Proof of Proposition~\ref{pro: Normalizable gci is hhci sometimes}}

Let $\nu:X \lra B$ be a polynomial normalization with Noether fibre
$F$. As before $F$ is connected and gci.

By assumption, the shearing map $\shr:HH^*(X|B) \lra HH^*(F)$ makes
$HH^*(F)$ into a finitely generated module over $HH^*(X|B)$. As we
saw in the proof of Theorem~\ref{thm: Finite gci spaces are hhci}
above, the shearing map $\psi:HH^*(F) \lra H_*(\Omega F)$ makes
$Z$ (the centre of $H_*(\Omega F)$) into a finitely generated
$HH^*(F)$-module. The relative shearing map
$\shr':HH^*(X|B) \lra H_*(\Omega F)$ is, by
Corollary~\ref{cor: relative shearing map is composition of two others}, the composition $\psi\shr$.
Thus $\shr'$ makes $Z$ into a finitely generated
$HH^*(X|B)$-module. By Theorem~\ref{thm: Sufficient condition},
$X$ is hhci.

To complete the proof of
Proposition~\ref{pro: Normalizable gci is hhci sometimes} we must
show that if $HH^*(X|B)$ is Noetherian then $\shr$ makes $HH^*(F)$
into a finitely generated $HH^*(X|B)$-module. This easily follows
from the next lemma.

\begin{lemma}
Set $Q=C^*(B)$, $R=C^*(X)$ and $S=C^*(F)$. Then $\Hom_{S^e}(S,S)$ is
a small $\Hom_{R^e_Q}(R,R)$-module.
\end{lemma}
\begin{proof}
Since $H^*(\Omega B)$ is finite, we see that $S$ is a small
$R$-module. As we saw earlier, the $R^e_Q$-module structure on
$S$ coming from the map of algebras $R^e_Q \lra R$ is the same as
the usual $R^e_Q$-module structure on $S$. 
Hence, $S$ is finitely built by $R$ over
$R^e_Q$ and therefore $\Hom_{R^e_Q}(R,S)$ is finitely built by
$\Hom_{R^e_Q}(R,R)$ as a $\Hom_{R^e_Q}(R,R)$-module.
\end{proof}

This completes the proof of
Proposition~\ref{pro: Normalizable gci is hhci sometimes}. \qed


\section{Some examples from group theory}
\label{sec:groupegs}

In this section we restrict attention to the special case
$R=C^*(BG;k)$, and seek to understand both the eci condition and the
gci condition.

\subsection{Chevalley groups.}
\label{subsec:Chevalley}
Now consider Quillen's $p$-adic construction of the classifying space
of a Chevalley group. As usual, all
spaces are completed at $p$, and we omit notation for this.
If  $q\neq p$ and $\Psi^q: BG \lra BG$ is an Adams map we
follow Quillen in defining
$F\Psi^q$ by the homotopy pullback square
$$\begin{array}{ccc}
F\Psi^q & \lra & BG \\
\downarrow && \downarrow B\Delta \\
BG &\stackrel{\{ 1, \Psi^q\} }\lra & BG \times BG
\end{array}$$
and then find a $p$-adic equivalence
$$F\Psi^q \simeq BG(q). $$
This gives a $p$-adic fibration
$$G \lra BG(q) \lra BG .$$

If $G$ is a sphere this shows that $BG(q) $ is an s-hypersurface space. 
For example, $G=U(1)$ this shows again 
that the classifying space $BGL_1(q)$ of the cyclic group $GL_1(q)$ is
an s-hypersurface space. More interesting
is the case $G=SU(2)$, which shows that  $BPSL_2(q)$ is an
s-hypersurface space.
When  $G$ is an iterated sphere bundle (e.g. if $G$ is one of the
classical groups) 
this shows  $BG(q)$ is wsci, and hence also gci and therefore eci.

\begin{problem}
Show that when $G=SU(3)$ and 
$q$ is such that the cohomology ring is not periodic,
then the space $BG(q)$  is not sci.
\end{problem}


Quillen \cite{QuillenGL}, Fiedorowicz-Priddy \cite{FP} and Kleinerman 
\cite{Kleinerman} show that $H^*(BG(q))$ is ci provided $H^*(BG)$ has
no $p$-torsion. Similarly Quillen shows the extraspecial groups 
\cite{extraspecial} have ci cohomology rings. As in
\ref{eg:homciisgci} this shows they are gci 
from the Eilenberg-Moore spectral
sequence and hence eci by Theorem \ref{thm:gciiseci}.


\subsection{Squeezed homology.}
Since we are working with groups, it is
illuminating to recall the first author's purely representation
theoretic calculation of the loop space homology $H_*(\Omega (BG_p^{\wedge}))$
\cite{squeezed}. In fact
$$H_*(\Omega(BG_p^{\wedge}))\cong \HOmega_*(G;k), $$
where $\HOmega_*(G;k) $ is defined algebraically.

More precisely $ \HOmega_*(G;k)$ is the homology of
$$\cdots \lra P_3 \lra P_2 \lra P_1 \lra P_0, $$
 a so-called {\em squeezed resolution}
of $k$. The sequence  of projective $kG$-modules $P_i$
is defined recursively as follows. To
start with $P_0=P(k)$ is the projective cover of $k$.
Now if $P_i$ has been constructed, take $N_i = \ker (P_i \lra P_{i-1})$
(where we take $P_{-1}=k$), and $M_i$ to be the smallest submodule of
$N_i$ so that $N_i/M_i$ is an iterated extension of copies of $k$.
Now take $P_{i+1}$ to be the projective cover of $M_i$.

\subsection{Some simple cases.}
Note that if $G$ is a $p$-group, we have $\Omega(BG_p^{\wedge})\simeq G$
so that the topology focuses on $H_*(\Omega BG) \cong kG$ and since
$k$ is the only simple module, $M_0=0$ and we again find $\HOmega_*(G)=kG$.

We would expect the next best behaviour to be when $H^*(BG)$ is a hypersurface.
Indeed, if $H^*(BG)$ is a polynomial ring modulo a relation of codegree $d$,
the Eilenberg-Moore spectral sequence
$$\Ext_{H^*(BG)}^{*,*}(k,k) \Rightarrow H_*(\Omega(BG_p^{\wedge}))$$
shows that there is an ultimate periodicity of period $d-2$. The actual
period therefore divides $d-2$. This same phenomenon can be seen in
the algebraic construction, where  complete information about
products is also available.

Here is an example where we can understand both the loop space
homology and the Hochschild cohomology explicitly.

\begin{example} ($G=p:q$ with $q|p-1$.) We take $G$ to be a
  non-trivial semidirect
product of $C_p$ with $C_q$, where $q$ is a divisor of $p-1$. Then
\[ H^*(BG,k)=\Lambda(y)\otimes k[x] \]
where $|y|=-(2q-1)$ and $|x|=-2q$ (we grade everything homologically). So
\[ H_*(\Omega BG\phat ;k) = \Lambda(\xi) \otimes k[\eta] \]
where $|\xi|=2q-1$ and $|\eta|=2q-2$. The spectral sequence
\cite{ShamirHHSS2} for the Hochschild  cohomology collapses to show
\[ HH^*(C^*(BG))=H^*(BG,H_*(\Omega BG\phat;k))=\Lambda(y,\xi)\otimes k[x,\eta]. \]
Alternatively, the spectral sequence
\[ HH^*(H^*(BG,k))\Rightarrow HH^*(C^*BG) \]
collapses, giving
\[ HH^*(C^*(BG))=\Lambda(y,\xi)\otimes k[x,\eta]. \]

Since the spectral sequence of~\cite{ShamirHHSS2} collapses on
the $E^2$-term, Theorem~\ref{thm: Sufficient condition} shows $BG$ is
hhci. Explicitly, we can take the element
$1\otimes \eta \in HH^*(BG)$ to show the hhci condition holds.
\end{example}

\begin{example} ($G=A_4$ with $p=2$.)
To start with we use a homotopy theoretic proof, showing that $BA_4$ is
an s-hypersurface space and hence also hhci.

Indeed, the natural 3-dimenional representation
$A_4 \lra SO(3)$ gives a 2-adic fibration
$$S^3\lra  BA_4 \lra BSO(3), $$
and $BA_4$ is an  s-hypersurface space at 2 with  $B\Gamma =BSO(3)$,
and $n=3$ (i.e., $d=-4$).

The stable cofibre sequence establishing that $BA_4$ is hhci will then be
$$BA_4 \times_{BSO(3)}BA_4 \lla BA_4 \lla \Sigma^2 BA_4,$$
and the periodicity element will be
$$\chi \in HH^{-2}(BA_4|BSO(3)).$$

Next we  outline a purely  algebraic proof.
To start with, we would like to see algebraically that
$H_*(\Omega (BA_4)_2^{\wedge})$ is eventually periodic. Although this calculation is already in \cite{squeezed}, we recall it briefly,
since we refer to it below.

This case is small enough to be able
to compute products in $H_*(\Omega BG\phat,k)$ using squeezed resolutions,
and we get
\[ H_*(\Omega BG\phat,k)=\Lambda(\alpha)\otimes
k\langle\beta,\gamma\rangle/(\beta^2,\gamma^2) \]
with $|\alpha|=1$ and $|\beta|=|\gamma|=2$.
Beware that $\beta$ and $\gamma$ do not commute, so that a $k$-basis
for $H_*(\Omega BG\phat,k)$ is given by alternating
words in $\beta$ and $\gamma$ (such as $\beta\gamma\beta$ or the empty word),
and $\alpha$ times these alternating words.

First note that
$$H^*(BA_4)=H^*(BV_4)^{A_4/V_4}=k[x_2,y_3,z_3]/(r_6)$$
where $r_6=x_2^3+y_3^2+ y_3z_3+z_3^2$.
From the Eilenberg-Moore spectral sequence we see that the loop space
homology will eventually have period dividing 4, and by calculation with
 the squeezed resolution we find the eventual period is exactly 4.

There are three simple modules. Indeed, the quotient of $A_4$ by its
normal Sylow 2-subgroup is of order 3; supposing for simplicity that $k$
contains three cube roots of unity $1, \omega, \omegabar$, the simples
correspond to how a chosen generator acts. The projective covers of the
three simple modules are
\[ P(k)=\vcenter{\xymatrix@=2mm{&k\ar@{-}[dl]\ar@{-}[dr]\\
\omega\ar@{-}[dr]&&\bar\omega\ar@{-}[dl]\\&k}},\quad
P(\omega)=\vcenter{\xymatrix@=2mm{&\omega\ar@{-}[dl]\ar@{-}[dr]\\
\bar\omega\ar@{-}[dr]&&k\ar@{-}[dl]\\&\omega}},\quad
P(\bar\omega)=\vcenter{\xymatrix@=2mm{&\bar\omega
\ar@{-}[dl]\ar@{-}[dr]\\k\ar@{-}[dr]&&\omega
\ar@{-}[dl]\\&\bar\omega}}. \]

Turning to Hochschild cohomology,
$HH^*(C^*(BG,k))$ is necessarily graded commutative, and we have
\[ HH^*(H^*(BG,k))=k[x,y,z]/(x^3+y^2+yz+z^2)\otimes
\Lambda(\xi,\eta,\zeta)\otimes k[\rho] \]
with $|\rho|=4$. We expect that both spectral sequences collapse.
\end{example}

\begin{example} ($G=L_3(2)$ at $p=2$.)
The behaviour is essentially the same as that of $G=A_4$ with $p=2$,
as described in \cite{squeezed}, and we give details here of a
different group by way of variation.

First
$$H^*(BL_3(2))=k [x_2,y_3, z_3]/(r_6),  $$
where $r_6=y_3z_3$.
From the Eilenberg-Moore spectral sequence we see that the loop space
homology will eventually have period dividing 4.

There are three simple modules, the trivial module $k$, and two others,
$M$ and $N$. The projective covers of the
three simple modules are
\[ P(k)=\vcenter{\xymatrix@=2mm{&k\ar@{-}[dl]\ar@{-}[dr]\\
M\ar@{-}[dr]&&N\ar@{-}[dl]\\&k}}, \quad
P(M)=\vcenter{\xymatrix@=2mm{&M\ar@{-}[dl]\ar@{-}[ddr]\\N\ar@{-}[d]\\
M\ar@{-}[d]&&k\ar@{-}[ddl]\\N\ar@{-}[dr]\\&M}},\quad
P(N)=\vcenter{\xymatrix@=2mm{&N\ar@{-}[dl]\ar@{-}[ddr]\\M\ar@{-}[d]\\
N\ar@{-}[d]&&k\ar@{-}[ddl]\\M\ar@{-}[dr]\\&N}}. \]
The squeezed resolution takes the following form, where the top row records
the projective modules $P_i$, and the second row records the modules $N_i$;
the modules $M_i$ are obtained by deleting the copies of $k$ marked with an
asterisk.

{\tiny
\[ \xymatrix@=1.2mm{\ar[r]&&\hspace{-4mm}P(M)\oplus P(N)\hspace{-4mm}
\ar[ddr] &\ar[rrr]&&&&\hspace{-4mm} P(N) \oplus P(M)\hspace{-4mm}
\ar[ddrr] &\ar[rr]&&&\hspace{-4mm} P(N)\oplus P(M)\hspace{-4mm}\ar[ddr]&
\ar[rr]&&&\hspace{-4mm}P(M)\oplus P(N)\hspace{-4mm} \ar[dr] &\ar[rr]&&&
\hspace{-4mm}P(k)\hspace{-4mm} \\ &&&&M\ar@{-}[d]&N\ar@{-}[d]
&&&&&&&N\ar@{-}[d]&&M\ar@{-}[d]&&M\ar@{-}[dr]&&N\ar@{-}[dl]\ar[ur] \\
k^*\ar@{-}[d]&k^*\ar@{-}[d]\ar[uur]&&&N\ar@{-}[d]&M\ar@{-}[d]&\ar[uur]&&&
k^*\ar@{-}[d]&k^*\ar@{-}[d]\ar[uur]&&\ M\ \ar@{-}[d]&&\ N\ \ar[uur]\ar@{-}[d]
&&&k& \\ M&N&&k^*\ar@{-}[dr]&M\ar@{-}[d]&N\ar@{-}[d]&k^*\ar@{-}[dl]&&&M&N&&
N\ar@{-}[d]&k^*\ar@{-}[dr]\ar@{-}[dl]&M\ar@{-}[d]\\ &&&&N&M&&&&&&&M&&N} \]}

Thus we see that
$$M_k \cong M_{k+2} \mbox{ for } k \geq 1$$

\end{example}

\begin{example}
\label{eg:BL33}
{\em (The cochains on $BL_3(3)$ at the prime $2$.)}
To start with, we note that this is also an example for $M_{11}$.
Indeed, the principal blocks of $L_3(3) $ and of $M_{11}$ are Morita
equivalent at the prime 2 \cite{ABG,Erdmann}, and there is a 2-adic  equivalence
$BL_3(3)\simeq BM_{11}$ \cite{Oliver}.

In any case, we take $X=BL_3(3)$ because we want to use the
Chevalley group properties.  This example is of interest since
$BL_3(3)$ is gci of codimension 1, whereas it is not wsci of
codimension 1. Since $BL_3(3)\simeq BSU_3(3)$
it is wsci of codimension 2 as in Subsection \ref{subsec:Chevalley}
in view of the 2-adic fibration
$$SU(3) \lra BL_3(3) \lra BSU(3). $$

Since it has cohomology ring
$$H^*(BL_3(3);k)\cong k[x_3,y_4,z_5]/(r_{10}), \mbox{ with }
r_{10}=x^2y+z^2, $$
we see that to be sci of codimension 1, we would need a space with
polynomial cohomology $k[x_3,y_4,z_5]$, which does not exist
\cite{AndersenGrodal}.

From the Eilenberg-Moore spectral sequence we see that the loop
space homology of $G$ will eventually have period dividing $8$, and we
may see this explicitly in terms of representation theory.

There are three simple modules: the trivial module $k$, a module $M$ of
dimension
$10$ and a module $N$ of dimension 44. The associated projective covers
$P(k), P(M)$ and $P(N)$ are as follows
\[ P(k)=\vcenter{\xymatrix@=2mm{&k\ar@{-}[dr]\ar@{-}[dl]\\
M\ar@{-}[d]\ar@{-}[ddrr]&&N\ar@{-}[d]\\k\ar@{-}[d]&&k\ar@{-}[d]\\
N\ar@{-}[dr]&&M\ar@{-}[dl]\\&k}},\quad
P(M)=\vcenter{\xymatrix@=2mm{&M\ar@{-}[dl]\ar@{-}[dr]\\
k\ar@{-}[d]\ar@{-}[ddrr]&&M\ar@{-}[d]\ar@{-}[ddll]\\N\ar@{-}[d]&&M\ar@{-}[d]\\
k\ar@{-}[dr]&&M\ar@{-}[dl]\\&M}},\quad
P(N)= \vcenter{\xymatrix@=2mm{&N\ar@{-}[dl]\ar@{-}[ddr]\\
k\ar@{-}[d]\\M\ar@{-}[d]&&N\ar@{-}[ddl]\\k\ar@{-}[dr]\\&N}},\quad \]

The squeezed resolution takes the following form, where the top row records
the projective modules $P_i$, and the second row records the modules $N_i$;
the modules $M_i$ are obtained by deleting the copies of $k$ marked with an
asterisk. The somewhat delicate part of the following calculation is the fact
that the $kk^*$ in $N_5$ has the effect of creating a module $M_5$ with the
same structure as $M_1$.

{\tiny
\[ \xymatrix@=1.2mm{\quad\ar[rr]&&&\hspace{-6mm}P(N)\oplus P(M)\hspace{-6mm}
\ar[ddr]&\ar[r]&&\hspace{-6mm}P(M)\oplus P(N)\hspace{-6mm}\ar[ddr]&\ar[rr]&&&
\hspace{-6mm}P(M)\oplus P(N)\hspace{-5mm}\ar[ddr]&\ar[r]&&\hspace{-5mm}P(N)
\oplus P(M)\hspace{-6mm}\ar[ddr]&\ar[rr]&&&\hspace{-6mm}P(M)\oplus P(N)
\hspace{-6mm}\ar[ddr]&\ar[rr]&&P(k)\\ &&M\ar@{-}[d]\ar@{-}[ddl]&&k^*\ar@{-}[d]
\ar@{-}[dr]&&&M\ar@{-}[d]\ar@{-}[ddr]&&&&k^*\ar@{-}[d]&&&&&M\ar@{-}[d]
\ar@{-}[ddl]&&M\ar@{-}[d]\ar@{-}[ddr]&N\ar@{-}[d]\\ &&M\ar@{-}[d]\ar[uur]&&
N\ar@{-}[d]&M\ar@{-}[d]\ar@{-}[ddl]\ar[uur]&&M\ar@{-}[d]&N\ar@{-}[d]\ar@{-}[dr]
&k^*\ar@{-}[d]\ar[uur]&&M\ar@{-}[d]&N\ar@{-}[ddl]\ar@{-}[d]\ar[uur]&&&&
M\ar@{-}[d]\ar[uur]&&k\ar@{-}[d]&k\ar@{-}[d]\ar[uur]\\ N\ar@{-}[d]&kk^*
\ar@{-}[dr]\ar@{-}[dl]&M\ar@{-}[d]&&k\ar@{-}[d]&k\ar@{-}[d]&&M\ar@{-}[d]&
k\ar@{-}[dl]&N&&k\ar@{-}[d]&k\ar@{-}[d]&&N\ar@{-}[d]&k\ar@{-}[dl]\ar@{-}[dr]&
M\ar@{-}[d]&&N\ar@{-}[d]&M\ar@{-}[dl]\\ N&&M&&M&N&&M&&&&N&M&&N&&M&&k\\} \]}
Thus we see
$$ M_k \cong M_{k+4} \mbox{ for } k \geq 1.$$
\end{example}


\end{document}